\documentclass[12pt,draft]{amsart}

\makeatletter

\theoremstyle{plain}    
\newtheorem{thm}{Theorem}[section]
\numberwithin{equation}{section} 
\numberwithin{figure}{section} 
\theoremstyle{plain}    
\newtheorem*{thm*}{Theorem} 
\theoremstyle{plain}    
\newtheorem{cor}[thm]{Corollary} 
\theoremstyle{plain}    
\newtheorem{lem}[thm]{Lemma} 
\theoremstyle{plain}    
\newtheorem{prop}[thm]{Proposition} 
\theoremstyle{definition}
\newtheorem{defn}[thm]{Definition}
\theoremstyle{remark}
\newtheorem{rem}[thm]{Remark}
\theoremstyle{remark}

\theoremstyle{remark}    

\theoremstyle{remark}    

\theoremstyle{definition}  
\newtheorem{example}[thm]{Example}
\theoremstyle{remark}
  \newtheorem*{acknowledgement*}{Acknowledgement} 

\theoremstyle{plain}    
\newtheorem{subthm}{Theorem}[subsection]
\theoremstyle{plain}    
\theoremstyle{plain}    
\newtheorem{sublem}[subthm]{Lemma} 
\theoremstyle{plain}    
\newtheorem{subprop}[subthm]{Proposition} 
\theoremstyle{definition}

\theoremstyle{remark}
\newtheorem{subrem}[subthm]{Remark}
\theoremstyle{remark}    

\theoremstyle{remark}    

\theoremstyle{plain}    

\newcommand{\id}{\operatorname{id}}

\newcommand{\TA}{\operatorname{T(A)}}
\newcommand{\UTAafd}{\operatorname{UT(A)_{AFD}}}
\newcommand{\UTAwafd}{\operatorname{UT(A)_{w-AFD}}}

\newcommand{\TAafd}{\operatorname{T(A)_{AFD}}}
\newcommand{\TAwafd}{\operatorname{T(A)_{w-AFD}}}

\usepackage{amssymb}
\usepackage{amscd}

\makeatother

\begin{document}

\title[tracial invariants]{tracial
invariants, classification and II$_1$ factor representations of Popa
algebras}

\author{Nathanial P. Brown}

\address{Department of Mathematics, Michigan State University, 
East Lansing, MI 48824}

\email{brown1np@cmich.edu}

\thanks{This research was supported by MSRI and NSF Postdoctoral 
Fellowships.  Currently on leave from Central Michigan University.}

\begin{abstract}
Using various finite dimensional approximation properties, four convex 
subsets of the tracial space of a unital C$^*$-algebra are defined.  One 
subset is characterized by Connes' hypertrace condition.  Another is 
characterized by hyperfiniteness of GNS representations.  The other two 
sets are more mysterious but are shown to be intimately related to 
Elliott's classification program. 

Applications of these tracial invariants include: 
\begin{enumerate}
\item An analogue of Szeg\"{o}'s Limit Theorem for arbitrary
self adjoint operators.

\item A McDuff factor embeds into $R^{\omega}$ if and only if it contains 
a weakly dense operator system which is injective.

\item There exists a simple, quasidiagonal, real rank zero
C$^*$-algebra with non-hyperfinite II$_1$ factor representations and
which is not tracially AF.  This answers negatively questions of Sorin
Popa and, respectively, Huaxin Lin.

\item If $A$ is any one of the standard examples of a stably
finite, non-quasidiagonal C$^*$-algebra and $B$ is a C$^*$-algebra
with Lance's WEP and at least one tracial state then there is  no
unital $*$-homomorphism $A \to B$. In particular, many stably finite,
exact C$^*$-algebras can't be embed into a stably finite, nuclear
C$^*$-algebra.
\end{enumerate}
\end{abstract}
\maketitle

\begin{center} 
{\em Dedicated to my big, beautiful family.}
\end{center}

\tableofcontents
\section{Introduction}

One of von Neumann's main motivations for initiating the study of
operator algebras was to provide a framework for studying unitary
representations of locally compact groups.  Hence it is no surprise
that representation theory of C$^*$-algebras has, historically,
received a lot of attention.  Over the last couple of decades,
however, representation theory has become less fashionable as exciting
new fields developed such as Connes' noncommutative geometry, Jones'
theory of subfactors, Elliott's classification program and, most
recently, Voiculescu's theory of free probability. 

In this paper we revisit representation theory but stick to the finite
(i.e.\ tracial) case.  There are two general questions which we study:
Approximation properties of traces on C$^*$-algebras and II$_1$ factor
representations of Popa algebras.  While we believe these topics to be
of independent interest, perhaps the most surprising part of this work
is that combining these two aspects of representation theory leads to
a variety of new results which provide a common thread between several
important problems in operator algebras.  Our results are most
strongly connected to questions in the classification program and free
probability but there are also relations with single operator theory
(more specifically, a classical theorem of Szeg\"{o} on spectral
distributions of self adjoint operators -- see Section 10), the
C$^*$-algebraic structure of the hyperfinite II$_1$ factor and even a
possible connection with the Novikov conjecture (though this is only
speculation at the moment).

Since we deal with some virtually disjoint subfields of operator
algebras, we feel it is worthwhile to paint a broad picture of the
topics covered before proving any results. 

\vspace{2mm}
\noindent{1.1. \bf Approximating Traces on C$^*$-algebras}

We will study four subsets of the tracial space of a unital,
separable C$^*$-algebra.  These subsets are defined via certain finite
dimensional approximation properties.  Though these approximation
properties may seem artificial at first glance, it turns out that they
are actually very natural.  Indeed, the definition of quasidiagonality
immediately leads to a norm approximation property for certain traces,
while the deep fact that the double dual of a nuclear C$^*$-algebra is
injective implies that {\em every} trace on a nuclear C$^*$-algebra
enjoys a strong 2-norm approximation property. 

One of the objectives of this paper is to point out how approximation
properties of traces carry information about (tracial) representation
theory of C$^*$-algebras.  This has already been witnessed in the work
of several authors (e.g.\ \cite{connes:classification},
\cite{connes:compactmetricspaces}, \cite{kirchberg:invent},
\cite{kirchberg:propertyTgroups}, \cite{bekka}, \cite{popa:simpleQD},
\cite{bedos:hypertraces}).  These papers are primarily concerned with
Connes' notion of a hypertrace.  For example, in \cite{bekka} Bekka
defines a unitary representation of a locally compact group to be
amenable if the C$^*$-algebra generated by the representation admits a
hypertrace.  (We recommend looking at \cite{bedos:hypertraces} for a
very nice introduction to the interactions between representation
theory and hypertraces.) It is a simple consequence of Voiculescu's
Theorem that any trace which satisfies one of the four approximation
properties defined in these notes will extend to a hypertrace in the
sense of Connes.  A wonderful result of Eberhard Kirchberg provides a
converse: If a trace extends to a hypertrace then it is `liftable'
(see \cite[Definition 3.1, Proposition
3.2]{kirchberg:propertyTgroups}).  It is a simple exercise to show
that liftable traces enjoy a weak 2-norm approximation property and
thus hypertraces are {\em characterized} by a certain finite
dimensional approximation property.

The three other approximation properties studied in this paper are 
natural variations on the approximation property enjoyed by hypertraces. 
One of these also turns out to be intimately related to representation 
theory: A trace is `uniformly weakly approximately finite dimensional' 
(see Definition \ref{thm:AFDtraces} below) if and only if the associated 
GNS representation gives a hyperfinite von Neumann algebra (see Theorem 
\ref{thm:mainthmUTAwafd}).  The remaining two approximation properties 
seem harder to understand, but it is natural to expect that they are 
also related to representation theory in a way which is not yet 
understood.  A better understanding of these other two sets of traces 
is closely related to questions about quasidiagonal C$^*$-algebras. 

The tracial invariants introduced here have a number of applications.
For example, in Section 10 we will see that approximation properties
of traces on type I C$^*$-algebras together with Voiculescu's Theorem
easily yield an analogue of a classical theorem of Szeg\"{o}.  We also
observe that they provide natural obstructions to the existence of
(unital) $*$-homomorphisms between certain classes of operator
algebras (cf.\ Corollary \ref{thm:nohomomorphisms}).  They also
clarify results of J.\ Rosenberg and S.\ Wassermann on the existence
of stably finite, non-quasidiagonal C$^*$-algebras (see the examples
at the end of Section 3 and the middle of section 7).

\vspace{2mm}
\noindent{1.2. \bf II$_1$ Factor Representations of Popa Algebras} 

Another goal of these notes is to study II$_1$ factor representations
of Popa algebras. 

\begin{defn}
A simple, separable, unital C$^*$-algebra, $A$, is called a {\em Popa
algebra} if for every finite subset $\mathfrak{F} \subset A$ and
$\varepsilon > 0$ there exists a nonzero finite dimensional
C$^*$-subalgebra $B \subset A$ with unit $e$ such that $\| ex - xe \|
\leq \varepsilon$ for all $x \in \mathfrak{F}$ and $e\mathfrak{F}e
\subset^{\varepsilon} B$ (i.e.\ for each $x \in \mathfrak{F}$ there
exists $b \in B$ such that $\| exe - b \| \leq \varepsilon$).
\end{defn}

Popa algebras are always quasidiagonal (QD). Using his local
quantization technique in the C$^*$-algebra setting Popa nearly
provides a converse in \cite{popa:simpleQD}: Every simple, unital,
quasidiagonal C$^*$-algebra with `sufficiently many projections'
(e.g.\ real rank zero) is a Popa algebra.  Thus the class of Popa
algebras is much larger than one might first guess. 

For some time there was speculation that quasidiagonality may be
closely related to nuclearity.  For example, in \cite[pg.\
157]{popa:simpleQD} Popa asked whether every Popa algebra with {\em
unique} trace is necessarily nuclear. (Counterexamples were first
constructed by Dadarlat in \cite{dadarlat:nonnuclearsubalgebras}.)
More generally, he asked in \cite[Remark 3.4.2]{popa:simpleQD} whether
the hyperfinite II$_1$ factor $R$ was the {\em only} II$_1$ factor
which could arise from a GNS representation of a Popa algebra.
Support for a positive answer was provided by the following
consequence of Connes' celebrated classification theorem:

\begin{thm*} 
Let $M$ be a (separable) II$_1$ factor.  Then $M \cong R$ if and only
if for every finite subset $\mathfrak{F} \subset M$ and $\varepsilon >
0$ there exists a nonzero finite dimensional C$^*$-subalgebra $B
\subset M$ with unit $e$ such that $\| ex - xe \|_2 \leq \varepsilon
\| e \|_2$ for all $x \in \mathfrak{F}$ and for each $x \in
\mathfrak{F}$ there exists $b \in B$ such that $\| exe - b \|_2 \leq
\varepsilon \|e\|_2$, where $\| \cdot \|_2$ is the 2-norm on $M$
coming from the unique trace.
\end{thm*}

Since Popa algebras are simple, the finite dimensional approximation
property from Definition 1.1 obviously passes to II$_1$ factor
representations (or any other representation).  Since the 2-norm
version of this approximation property characterizes $R$ it was
natural to expect that the II$_1$ factor representation theory of Popa
algebras would be trivial.  We will see, however, that the
representation theory of Popa algebras is actually very rich.  For
example, we will construct a Popa algebra $A$ with the property that
for every (separable) II$_1$ factor $M$ there exists a tracial state
$\tau$ on $A$ such that $\pi_{\tau} (A)^{\prime\prime} \cong M
\bar{\otimes} R$ (see Theorem \ref{thm:arbitraryMcDuff}).  On the
other hand we will see that if $A$ is a locally reflexive (e.g.\
exact) Popa algebra with {\em unique} trace $\tau$ then $\pi_{\tau}
(A)^{\prime\prime} \cong R$ thus giving a positive answer to Popa's
question in this case (see Corollary
\ref{thm:locallyreflexiveuniquetrace}).

\vspace{2mm}
\noindent{1.3. \bf Applications} 

As mentioned above, our study of finite representation theory leads to
a number of results which are closely related to some important open
problems. 

\vspace{2mm}
{\noindent\bf Elliott's Classification Program} 

There are three new results which those in the classification program
may find of interest.  First, we will show that if Elliott's
conjecture holds for simple, nuclear, QD C$^*$-algebras then every
trace on every nuclear, QD (not necessarily simple!) C$^*$-algebra
satisfies the strongest kind of approximation property (cf.\
Proposition \ref{thm:elliott}).  Hence we get a new necessary
condition for Elliott's conjecture to hold.  Second, we will point out
why this necessary condition may also eventually become part of a
sufficient condition for classification (see Theorem
\ref{thm:idontknow}).

The third result of relevance to the classification program is
inspired by Huaxin Lin's class of tracially AF algebras.  Roughly
speaking a tracially AF algebra is a Popa algebra with the additional
property that the finite dimensional algebra $B$ from Definition 1.1
can always be taken `large in trace' (see \cite{lin:TAF} for the
precise definition).

Thanks to Huaxin Lin's remarkable classification result for tracially AF
algebras (cf.\ \cite{lin:TAFclassification}) it is now of fundamental
importance to understand when a Popa algebra is a tracially AF algebra.
Indeed, Elliott's classification conjecture predicts that every {\em
nuclear} Popa algebra with real rank zero, unperforated K-theory and
Riesz decomposition property must be tracially AF (see Proposition
\ref{thm:elliottpredictsTAF}). Thus, by Lin's classification theorem,
verifying this case of Elliott's conjecture is {\em equivalent} to the
following question (modulo a UCT assumption):

\vspace{2mm}
{\em Does every nuclear Popa algebra with real rank zero, unperforated 
K-theory and having the Riesz decomposition property also have 
tracial topological rank zero, in the sense of Huaxin Lin
\cite{lin:tracialtopologicalrank}?}
\vspace{2mm}

Using the tracial invariants alluded to above and representation
theory of Popa algebras we will show that it is possible to construct
an {\em exact} Popa algebra with virtually every nice property
(including all those above) but which is {\em not} of tracial
topological rank zero (cf.\ Corollary \ref{thm:notTAF}).  This example
is a bit surprising as some experts did not expect that the question
above would have much to do with nuclearity.  (See, for example,
\cite[page 694]{lin:TAF} where it was ``tempting to conjecture that
every quasidiagonal, simple C$^*$-algebra of real rank zero, stable
rank one and with weakly unperforated $K_0$ is tracially AF.'')
Indeed, the question above can be regarded as a finite analogue of the
question of whether every {\em nuclear}, simple, infinite
C$^*$-algebra is purely infinite and hence our example in Corollary
\ref{thm:notTAF} should be regarded as a finite analogue of M.\
R{\o}rdam's recent construction of a (non-nuclear) simple, infinite
but not purely infinite C$^*$-algebra (see
\cite{rordam:infinitenotpurelyinfinite}).

\vspace{2mm}
{\noindent\bf Free Probability} 

There are two new results arising from this work which are relevant to
free probability.  One is related to Connes' embedding problem and the
other is related to the semicontinuity question for Voiculescu's free
entropy dimension.  

As mentioned above, we will show that every McDuff factor contains a
weakly dense Popa algebra.  It seems quite likely that many other
II$_1$ factors contain weakly dense Popa algebras as well and hence it
is natural to ask how free entropy dimension behaves on Popa algebras.
For example, is it true that for Popa algebras the free entropy
dimension with respect to {\em any} trace and any set of generators is
{\em always} $\leq 1$?  The point is that Popa algebras do not arise
from any sort of (reduced) free product construction (since such free
products are essentially never QD) and have abstract properties which
are quite similar to properties which characterize the hyperfinite
II$_1$ factor.  Thus it is important to determine how free entropy
dimension behaves on these algebras.

Regarding Connes' embedding problem (i.e.\ the question of whether or
not microstates always exist) we obtain the following result (see
Theorem \ref{thm:embeddableMcDuff}): A II$_1$ factor $M$ has
microstates if and only if there exists a weakly dense operator system
$X \subset R \bar{\otimes} M$ which is injective.  Hence we see that
the difference between embedding into the ultrapower of $R$ and
actually being isomorphic to $R$ is rather delicate.  Our own feeling
is that this difference is too delicate to expect that every II$_1$
factor has microstates, but we have not yet been able to construct a
counterexample.

\vspace{2mm}
{\noindent\bf Szeg\"{o}'s Limit Theorem}

Inspired by ideas of Arveson (see also B\'{e}dos' work in this
direction; \cite{bedos:Szego}) we observe in Section 10 that the
techniques and ideas of this paper easily yield a very general
existence theorem, analogous to the classical limit theorem of
Szeg\"{o}, for {\em arbitrary} self adjoint operators on a separable
Hilbert space.  In particular, Theorem \ref{thm:Szego} is a vast
generalization of \cite[Theorem 4.5]{arveson:numerical} in the sense
that {\em all} of the hypotheses regarding the larger C$^*$-algebra
are unnecessary.

\vspace{2mm}
{\noindent\bf Quasidiagonality and the Hyperfinite II$_1$ Factor}

It is an open question whether or not $R$ is a quasidiagonal (QD)
C$^*$-algebra.  While many experts believe that $R$ is not QD, there
is little concrete evidence to support a negative answer.  On the
other hand, a positive answer would imply, for example, that every
simple, unital, nuclear, stably finite C$^*$-algebra is QD (something
predicted by Elliott's conjecture) and that $C_r^* (\Gamma)$ is QD for
every discrete amenable group $\Gamma$ (a conjecture of J.\
Rosenberg).  Hence, aside from being a natural and basic problem, we
feel this to be an important question as well.  In the appendix we
will reformulate this problem in terms of approximation properties of
traces on (separable) C$^*$-algebras.  Thus it is our hope that future
work on approximation properties of traces will shed light on this
problem.

\vspace{2mm}

The paper is organized as follows.  In section 2 we collect a number
of preliminary results which we will need.  The main result of this
section shows how to construct Popa algebras with specified GNS
representations.  In section 3 we define various subsets of traces and
study their properties.  This section also contains a number of
examples.  In section 4 we study the case that $A$ is a nuclear
C$^*$-algebra.  It is shown  that Elliott's conjecture
predicts that every trace on every nuclear, quasidiagonal
C$^*$-algebra (simple or not) has the strongest kind of approximation
property.  In section 5 we consider the exact case and show how to
construct exact Popa algebras which are not tracially AF.  Section 6
treats the locally reflexive case while Section 7 treats the WEP
case. In section 8 we discuss II$_1$ factor representations of Popa
algebras and observe that even very nice Popa algebras can have
non-hyperfinite II$_1$ factor representations. Section 9 
observes that the techniques of this paper easily yield new
characterizations of McDuff factors which are embeddable into the
ultrapower of the hyperfinite II$_1$ factor.  One consequence of this 
result is that a number of well known II$_1$ factors (e.g.\ $R\bar{\otimes} 
L(\Gamma)$ for any i.c.c., residually finite, discrete group) contain weakly 
dense injective subspaces, but, of course, are not themselves 
injective.  Section 10 contains an analogue of Szeg\"{o}'s Limit 
Theorem.  Section 11 just contains a number of questions which arise 
naturally from this work.  Finally, we have an appendix which discusses 
the question of whether or not the hyperfinite II$_1$ factor is 
a quasidiagonal C$^*$-algebra.

\begin{acknowledgement*}
This work began during the year-long program in
operator algebras at MSRI, 2000-2001.  I spoke to nearly everyone I
encountered during that year about various aspects of this work and it
would be impossible to recall all of the people who contributed
remarks and ideas.  Instead I express my sincerest thanks to MSRI,
the organizers and all of the participants for producing such an open
and stimulating research environment.  However, I do need to
specifically thank Marius Dadarlat and Dimitri Shlyakhtenko for a
(seemingly infinite) number of helpful discussions.
\end{acknowledgement*}

\section{Preliminaries}
In this section we present a number of  results (many of which are 
well known, but restated in particularly convenient ways) which serve
as the backbone for the rest of the paper. The main new result, Theorem
\ref{thm:basicconstruction}, states that residually finite dimensional
C$^*$-algebras always admit approximately trace preserving embeddings
into Popa algebras.

\subsection{Notation}
Before stating any results we wish to introduce some
notation which will be used throughout this paper.  

{\em We first remark that, unless otherwise noted or obviously 
false, all C$^*$-algebras in 
this paper are assumed to be unital and separable.}  Similarly, all 
von Neumann algebras will be assumed to have separable preduals (with 
the exception of $R^{\omega}$, which is well known to be non-separable).

For a Hilbert
space $H$, we will let $B(H)$ and $\mathcal{K}(H)$ denote the bounded
and, respectively, compact operators on $H$.  We let $\| \cdot \|$
denote the operator norm on $B(H)$ while $\| \cdot \|_{HS}$ and
$<\cdot ,\cdot >_{HS}$ will denote the Hilbert-Schmidt norm and,
respectively, inner product on the Hilbert-Schmidt operators on $H$.

When $A$ is a C$^*$-algebra with state $\eta$ we will denote the
associated GNS Hilbert space, representation and von Neumann algebra 
 by $L^2 (A,\eta)$, $\pi_{\eta} : A \to B(L^2 (A,\eta))$ and 
$\pi_{\eta}(A)^{\prime\prime}$, respectively.  

The symbols $\odot$, $\otimes$ and $\bar{\otimes}$ will denote the
algebraic, minimal and W$^*$-tensor products, respectively.

If $A$ is a C$^*$-algebra we will
let $A^{op}$ denote the opposite algebra (i.e.\ $A^{op} = A$ as
involutive normed linear spaces, but multiplication in $A^{op}$ is
defined by $a \circ b = ba$; the latter multiplication being the given
multiplication in $A$).  $A^{**}$ will denote the enveloping von
Neumann algebra of $A$ (i.e.\ the Banach space double dual of $A$).
Contrary to our standing assumption that von Neumann algebras should
have separable preduals, $A^{**}$ often has a non-separable predual.

If $\tau \in \TA$ is a tracial state then there is a canonical
antilinear isometry $J: L^2(A,\tau) \to L^2(A,\tau)$ defined by
$J(\hat{a}) = \hat{a^*}$.  One defines a $*$-homomorphism
$\pi_{\tau}^{op} : A^{op} \to B(L^2(A,\tau))$ by $\pi_{\tau}^{op} (a)
= J\pi_{\tau} (a^*)J$. Since $J\pi_{\tau}(A)J \subset
\pi_{\tau}(A)^{\prime}$ one then gets an algebraic homomorphism
$\pi_{\tau} \odot \pi_{\tau}^{op} : A \odot A^{op} \to B(L^2(A,\tau))$
defined on elementary tensors by $\pi_{\tau} \odot \pi_{\tau}^{op}
(a\otimes b) = \pi_{\tau} (a) \pi_{\tau}^{op} (b)$.

Completely positive
maps (cf.\ \cite{paulsen:cbmaps}) will play an important role in these
notes.  We will use the abbreviations c.p. and u.c.p. for `completely
positive' and `unital completely positive', respectively.  

The hyperfinite II$_1$ factor will appear many times and will always
be denoted by $R$.  We will use $tr_n$ to denote the unique tracial
state on the $n\times n$ matrices.  For a von Neumann algebra, $M$,
with faithful, normal, tracial state $\tau$, we will let $\| \cdot
\|_{2,\tau}$ be the associated 2-norm (i.e.\ $\| x \|_{2,\tau} =
\tau(x^* x)^{1/2}$).  If $M$ has a unique trace (i.e.\ is a factor)
then we will drop the dependence on $\tau$ and simply write $\| \cdot
\|_2$.

Finally, we will need the ultrapower of the hyperfinite II$_1$ factor.
That is, given a free ultrafilter $\omega \in \beta {\mathbb N}
\backslash {\mathbb N}$ one defines an ideal $I_{\omega} \subset
l^{\infty}(R) = \{ (x_n) \in \Pi_{n \in {\mathbb N}} R: \sup_{n \in
{\mathbb N}} \| x_n \| < \infty \}$ by $I_{\omega} = \{ (x_n) \in
l^{\infty}(R) : \lim_{n \to \omega} \|x_n \|_2 = 0 \}$.  Then the
ultrapower of $R$ with respect to $\omega$ is defined to be the
(C$^*$-algebraic) quotient: $R^{\omega} = l^{\infty}(R)/I_{\omega}$.
$R^{\omega}$ is a II$_1$ factor with trace $\tau_{\omega} ((x_n) +
I_{\omega}) = \lim_{n \to \omega} \tau_{R} (x_n)$.

\subsection{Voiculescu's Theorem}

We will need the following version of Voiculescu's Theorem.  

\begin{subthm}
\label{thm:Voiculescu'sthm}
Let $A \subset B(H)$ be in general position (i.e.\ $A \cap
\mathcal{K}(H) = \{ 0 \}$).  If $\phi : A \to M_n ({\mathbb C})$ is a
u.c.p. map then there exist isometries $V_k : {\mathbb C}^n \to H$
such that $\| \phi(a) - V_k^* a V_k \| \to 0$, for all $a \in A$, as
$k \to \infty$.  Moreover if $P_k = V_k V_k^*$ then one has the
following estimates on commutators:

\begin{enumerate} 
\item $\limsup \| P_k a - a P_k\| \leq \max \{ \|
\phi(aa^*) - \phi(a)\phi(a^*) \|^{1/2}, \| \phi(a^*a) -
\phi(a^*)\phi(a) \|^{1/2} \}$,

\item $\limsup \frac{\| P_k a - a P_k\|_{HS}}{\| P_k \|_{HS}} \leq
\bigg(  tr_n \big( \phi(aa^*) - \phi(a)\phi(a^*)\big)  +  tr_n
\big(\phi(a^*a) - \phi(a^*)\phi(a)\big)  \bigg)^{1/2}$.
\end{enumerate}
\end{subthm} 
 
\begin{proof}  
This result is really contained in the standard proof of Voiculescu's
Theorem (cf.\ \cite{arveson:extensions} or \cite{davidson}).  However,
the commutator estimates, which are the important part for us, become
especially transparent if we use the usual versions of Voiculescu's
Theorem and Stinespring's Theorem.

So, let $\pi : A \to B(K)$ be the Stinespring dilation of $\phi$ with
isometry $V : {\mathbb C}^n \to K$ such that $\phi(a) = V^* \pi(a)
V$. Define the projection $P = VV^*$ and note that we have the
identity $P\pi(a) - \pi(a)P = Pa(1 -P) - (1 - P)aP$.  Note also that
$Pa(1 -P)$ and $(1 - P)aP$ have orthogonal domains and ranges and, in
particular, are perpendicular in the Hilbert space of Hilbert-Schmidt
operators.  Using this one can verify the identities $$\| P \pi(a) -
\pi(a) P\| = \max \{ \| \phi(aa^*) - \phi(a)\phi(a^*) \|^{1/2}, \|
\phi(a^*a) - \phi(a^*)\phi(a) \|^{1/2} \},$$ and $$\frac{\| P \pi(a) -
\pi(a) P \|_{HS}}{\| P \|_{HS}} = \bigg( tr_n \big( \phi(aa^*) -
\phi(a)\phi(a^*)\big) + tr_n \big(\phi(a^*a) - \phi(a^*)\phi(a)\big)
\bigg)^{1/2}.$$

Note that the same estimates hold for the representation $\iota \oplus
\pi : A \to B(H\oplus K)$, where $\iota : A \hookrightarrow B(H)$ is
the given inclusion.  Since, $\iota$ is approximately unitarily
equivalent to $\iota \oplus \pi$, it is clear how to complete the
proof.
\end{proof}

In one place, we will need the following technical version of Voiculescu's 
Theorem.  A proof can be found in \cite{brown:QDsurvey} or 
\cite{brown:herrero}. 

\begin{subprop}
\label{thm:technicalVoiculescuThm}
Let $A \subset B(H)$ be in general position and $\Phi : A \to B(K)$ 
be a u.c.p. map which is a faithful $*$-homomorphism modulo the 
compacts (i.e.\ composing with the quotient map to the Calkin algebra 
yields a faithful $*$-monomorphism $A \hookrightarrow Q(K)$).  Then 
there exists a sequence of unitaries $U_n : K \to H$ such that for 
every $a \in A$ we have  
$$\limsup \| a - U_n \Phi(a) U_n^* \| \leq 2\max \{ 
\| \Phi(aa^*) - \Phi(a)\Phi(a^*) \|^{1/2} , 
\| \Phi(a^* a) - \Phi(a^*)\Phi(a) \|^{1/2} \}.$$ 
\end{subprop}

\subsection{Elliott's Intertwining Argument} 

Perhaps the single most important argument in the classification
program is due to George Elliott.  Though usually done in the setting
of C$^*$-algebras we will need Elliott's approximate intertwining
argument in the setting of von Neumann algebras.  While not the most 
general possible form, the following version is more than sufficient 
for our purposes.  The set-up is as follows. 

Assume that $M \subset B(L^2 (M,\tau))$ and $N \subset
B(L^2(N,\gamma))$ are von Neumann algebras acting standardly, with
faithful, normal tracial states $\tau$ and, respectively, $\gamma$.
Let $X_1 \subset X_2 \subset \ldots \subset M$ and $Y_1 \subset Y_2
\subset \ldots \subset N$ be (not necessarily unital)
C$^*$-subalgebras such that $\cup X_i$ is weakly dense in $M$ and
$\cup Y_i$ is weakly dense in N.  Further assume that we have
c.p. maps $\alpha_n : Y_n \to X_n$, $\beta_n : X_n \to Y_{n + 1}$,
which are contractive both with respect to the operator norms and the
2-norms coming from $\tau$ and $\gamma$, and finite
subsets $\Lambda_i \subset X_i$ and $\Omega_i \subset Y_i$ with the
following properties:

\begin{enumerate}
\item $\Lambda_i \subset \Lambda_{i + 1}$, $\Omega_i \subset \Omega_{i
+ 1}$, for all $i \in {\mathbb N}$, and the linear spans of $\cup
\Lambda_i$ and $\cup \Omega_i$ are norm dense in $\cup X_i$ and $\cup
Y_i$, respectively, and hence weakly dense in $M$ and $N$,
respectively.  To simplify things, we will also assume that $x_1, x_2
\in \Lambda_i \Longrightarrow x_1x_2 \in \Lambda_{i+1}$ and,
similarly, that $\Omega_{i+1}$ contains the product of any pair of
elements from $\Omega_i$.

\item $\alpha_i (\Omega_i) \subset \Lambda_i$ and $\beta_i (\Lambda_i) 
\subset \Omega_{i + 1}$ for all $i \in {\mathbb N}$.  

\item Both $\{ \alpha_i \}$ and $\{ \beta_i \}$ are weakly
asymptotically multiplicative. That is, $\| \alpha_i (y_1 y_2) -
\alpha_i (y_1) \alpha_i (y_2) \|_{2,\tau} \to 0$, as $i \to \infty$,
for all $y_1, y_2 \in \cup Y_i$ and similarly for $\{ \beta_i \}$.
\end{enumerate}

\begin{subthm}(Elliott's Intertwining)
In the setting described above, if it happens that $\| x - 
\alpha_{n + 1} \circ \beta_n (x) \|_{2,\tau} < 1/2^n$ and $\| y - 
\beta_n \circ \alpha_n (y) \|_{2,\gamma} < 1/2^n$ for all $x \in 
\Lambda_n$, $y \in \Omega_n$ and all $n \in {\mathbb N}$, then 
$M \cong N$.
\end{subthm}

\begin{proof} As this argument is well known, we will be a bit 
sketchy.  Two facts which we will need are: $i)$ every norm bounded
sequence which is Cauchy in 2-norm converges (in 2-norm) and $ii)$ on
norm bounded subsets, the 2-norm topology is the same as the strong
operator topology (since our von Neumann algebras are acting
standardly on the $L^2$ spaces coming from their traces).

 Since the $\alpha_n$'s and $\beta_n$'s are 2-norm 
contractive, one first checks that for each $y \in \cup
\Omega_i$, the sequence $\{ \alpha_n (y) \}$ is Cauchy in 2-norm 
(and similarly for each $x \in \cup \Lambda_i$).  Hence this
is also true on the linear spans of $\cup \Lambda_i$ and $\cup
\Omega_i$.

Since the $\alpha_n$'s and $\beta_n$'s are norm contractive, it
follows that for each $y$ in the linear span of $\cup \Omega_i$, the
sequence $\{ \alpha_n (y) \}$ is convergent in $M$ (and similarly for
all $x \in span(\cup \Lambda_i$)).  Hence we can define linear maps
$\Phi : span(\cup \Lambda_i) \to N$ and $\Psi : span(\cup \Omega_i)
\to M$ by $\Phi(x) = \lim \beta_n(x)$ and $\Psi (y) = \lim
\alpha_n(y)$.  Note that $\Phi$ and $\Psi$ are contractive with
respect to both the operator norms and the 2-norms.  This implies that
$\Phi$ and $\Psi$ can be (uniquely) extended to the norm closures of
$span(\cup \Lambda_i)$ and $span(\cup \Omega_i)$ (which are weakly
dense C$^*$-algebras, by condition (1) above) and, moreover, that
these extensions are 2-norm contractive as well.  Now one uses
Kaplansky's density theorem (and the fact that our extensions are
still norm and 2-norm contactive) to extend beyond these weakly dense
C$^*$-subalgebras to all of $M$ and $N$. (i.e.\ For each $x \in M$ we
take a norm bounded sequence $\{ x_n \}$ from the norm closure of
$span(\cup \Lambda_i)$ which converges to $x$ in 2-norm. The image of
$\{ x_n \}$ in $N$ is then a norm bounded sequence which is Cauchy in
2-norm and hence we map $x$ to the (2-norm) limit of this sequence.)

To save notation, we will also let $\Phi : M \to N$ and $\Psi : N \to
M$ denote the maps constructed in the previous paragraph.  Note that
these maps are 2-norm contractive and linear.  They are also
$*$-preserving since the strong and strong-$*$ topologies agree on
bounded subsets of a tracial von Neumann algebra (since $\| x
\|_{2,\tau} = \| x^* \|_{2,\tau}$).  It is also easy to check that they
are mutual inverses on the spans of $\cup \Lambda_i$ and $\cup
\Omega_i$.  By 2-norm contractivity it follows that they are mutual
inverses on all of $M$ and $N$.  Hence we only have to observe that
both $\Phi$ and $\Psi$ are multiplicative on $M$ and $N$.  Since
multiplication is continuous, on bounded sets, in the 2-norm, it is
not hard to check that $\Phi$ and, respectively, $\Psi$ are
multiplicative on the norm closures of the linear spans of $\cup
\Lambda_i$ and, respectively, $\cup \Omega_i$.  Finally, another
application of Kaplansky's density theorem and a standard
interpolation argument allow one to deduce multiplicativity on all of
$M$ and $N$.
\end{proof}

\subsection{Consequences of Elliott's Conjecture}
We remind the reader that all algebras are assumed separable and unital. 

We state here the special case of Elliott's Conjecture which will be
relevant for us.  We then deduce a few statements which are predicted
by this conjecture.  We are indebted to M. R{\o}rdam and H. Lin for
some helpful discussions regarding these issues.

For a stably finite C$^*$-algebra $A$, the Elliott invariant is the
triple $(K_0 (A), K_1 (A), \TA)$, where $\TA$ is the set of tracial
states on $A$, together with the natural pairing $P_A : K_0 (A) \times
\TA \to {\mathbb R}$.  Given two algebras $A$ and $B$, we say that
their Elliott invariants are isomorphic if $K_1 (A) \cong K_1(B)$ and
there exist a scaled, ordered group isomorphism $\Phi : K_0 (A) \to
K_0 (B)$ and an affine homeomorphism $T : \TA \to T(B)$ such that
$P_A(x,\tau) = P_B (\Phi(x),T(\tau))$, for all $(x,\tau) \in K_0 (A)
\times \TA$.

(Special case of) {\bf Elliott's Conjecture:} If two simple, stably
 finite, nuclear C$^*$-algebras have isomorphic Elliott invariants (as
 described above) then they are isomorphic.

\begin{subprop}
\label{thm:elliottpredictsASH}
Elliott's conjecture predicts that if $A$ is a stably finite, simple,
nuclear C$^*$-algebra then for any UHF algebra, $\mathcal{U}$, $A
\otimes \mathcal{U}$ is an inductive limit of subhomogeneous algebras
(i.e.\ ASH).
\end{subprop}

\begin{proof} 
By \cite[Theorem 5.2 (b)]{rordam:tensorUHFII} it follows that the
pairing $P_{A \otimes \mathcal{U}} : K_0(A \otimes \mathcal{U}) \times
T(A \otimes \mathcal{U}) \to {\mathbb R}$ is weakly unperforated
(which means that if $P_{A \otimes \mathcal{U}} (x,\tau) > 0$ for all
$\tau$ then $x > 0$).  Note that in order to apply R{\o}rdam's results
from \cite{rordam:tensorUHFII} we have to rely on U. Haagerup's result
that quasitraces on (unital) exact C$^*$-algebras are traces (cf.\
\cite{haagerup:quasitraces}).  Once we know that the invariant of $A
\otimes \mathcal{U}$ is weakly unperforated, we are done since Elliott
showed how to construct a simple ASH algebra with arbitrary weakly
unperforated invariant (see, for example, the appendix of
\cite{elliott-villadsen}).  Indeed, if Elliot's conjecture holds, we
can find an ASH algebra whose invariant is isomorphic to the invariant
of $A \otimes \mathcal{U}$ and hence $A \otimes \mathcal{U}$ is
isomorphic to an ASH algebra.
\end{proof}

Note that Elliott's conjecture also predicts that every simple, stably
finite, nuclear C$^*$-algebra is QD, since quasidiagonality passes to
subalgebras and ASH algebras are QD.

\begin{subprop}
\label{thm:elliottpredictsTAF}
Elliott's conjecture predicts that if $A$ is simple, stably finite,
nuclear, real rank zero (cf.\ \cite{brown-pedersen}), has weakly
unperforated invariant (see the proof of the previous proposition for
this definition) and $K_0(A)$ has the Riesz interpolation property
(cf.\ \cite[Section IV.6]{davidson}) then $A$ is tracially AF in the
sense of \cite{lin:TAF}.
\end{subprop}

\begin{proof}  We first remark that in the real rank zero case the tracial 
simplex is no longer relevant and hence Elliott's invariant reduces to
K-theory alone.  Indeed, if both $A$ and $B$ are C$^*$-algebras of
real rank zero and $\Phi : K_0 (A) \to K_0 (B)$ is a scaled, ordered
group isomorphism such that $\Phi([1_A]) = [1_B]$ then $\Phi$ induces
an affine homeomorphism $\TA \to T(B)$ since we may (affinely,
homeomorphically) identify $\TA$ (resp.\ $T(B)$) with the states in
${\rm Hom}(K_0(A), {\mathbb R})$ (resp.\ ${\rm Hom}(K_0(B), {\mathbb
R})$). To see that this is true, we first note that the obvious map
$\TA \to {\rm Hom}(K_0(A), {\mathbb R})$ is affine and injective since
$A$ has real rank zero.  It is also onto the states in ${\rm
Hom}(K_0(A), {\mathbb R})$ since every state on $K_0 (A)$ comes from a
trace on $A$ when $A$ is unital and exact (cf.\
\cite{haagerup-thorbjornsen}).  Finally, it is easy to check (again
using real rank zero) that a sequence of traces $\tau_n \in \TA$
converges to $\tau \in \TA$ in the weak-$*$ topology if and only if
their images in ${\rm Hom}(K_0(A), {\mathbb R})$ converge in the
topology of pointwise convergence and hence our identification is also
a homeomorphism.

Finally, in \cite{elliott-gong} it is shown how to construct simple AH
algebras with real rank zero and with arbitrary unperforated K-theory
and Riesz interpolation property.  As observed by Lin,
\cite[Proposition 2.6]{lin:TAF}, the Elliott-Gong construction always
yields tracially AF algebras and hence we can find a simple, nuclear,
tracially AF algebra with the same K-theory as $A$.  Hence, if
Elliott's conjecture holds, $A$ is isomorphic to a tracially AF
algebra.
\end{proof}

\subsection{Approximately Trace Preserving Embeddings into Popa Algebras}

We now present the main technical result of this paper.  The informed
reader will note that every aspect of this result can be traced back
to the classification program.  Indeed, we will adapt the inductive
limit techniques of Dadarlat \cite{dadarlat:nonnuclearsubalgebras} to
construct new Popa algebras and use Elliott's intertwining argument to
understand their GNS representations.

In the following theorem $\mathfrak{C}$ will denote some collection of
C$^*$-algebras which is closed under increasing unions (i.e.\
inductive limits with injective connecting maps) and tensoring with
finite dimensional matrix algebras.  A C$^*$-algebra $E$ is called
{\em residually finite dimensional} if $E$ has a separating family of
finite dimensional representations (i.e.\ for every $0 \neq x \in E$
there exists a $*$-homomorphism $\pi : A\to M_n ({\mathbb C})$ such
that $\pi(x) \neq 0$).

\begin{subthm}
\label{thm:basicconstruction}
Let $E \in \mathfrak{C}$ be a residually finite dimensional
C$^*$-algebra.  Then there exists a Popa algebra $A \in \mathfrak{C}$
such that for every $\varepsilon > 0$ we can find a $*$-monomorphism
$\rho : E \hookrightarrow A$ with the property that for each trace
$\tau \in T(E)$, there exists a trace $\gamma \in \TA$ such that 

\begin{enumerate}
\item $|\gamma\circ\rho (x) - \tau(x) | < \varepsilon \| x \|$ for
all $x \in E$ and,

\item $\pi_{\gamma} (A)^{\prime\prime} \cong R \bar{\otimes} \pi_{\tau}
(E)^{\prime\prime}.$
\end{enumerate}
\end{subthm}

The proof of this result becomes much more transparent once the main
idea is understood.  Hence we think it is worthwhile to give the main
idea first and leave the details to the end.

So suppose that $E$ is a residually finite dimensional C$^*$-algebra
and $\tau \in T(E)$. Let $\mathcal{U}$ be some UHF algebra. Then the
canonical, unital inclusion $E \hookrightarrow E\otimes \mathcal{U}$
is honestly trace preserving (in fact, yields an isomorphism of
tracial spaces) and the weak closure in any GNS representation is
obviously of the form $R \bar{\otimes} \pi_{\tau} (E)^{\prime\prime}.$
The problem, of course, is that $E\otimes \mathcal{U}$ is not a Popa
algebra.  So the idea is that we will use an inductive limit
construction to get a sequence $$ E \to E\otimes M_{k(1)} ({\mathbb
C}) \to E\otimes M_{k(1)} ({\mathbb C}) \otimes M_{k(2)} ({\mathbb C})
\to \cdots,$$ such that the limit is a Popa algebra, but the
connecting maps above will be chosen so that when one applies a trace
it will (approximately) look like the sequence which yields $E \otimes
\mathcal{U}$.

We now describe the basic construction which will be needed to get our
Popa algebras.  Let $\pi : E \to M_k ({\mathbb C})$ be a
representation, $\tau \in T(E)$ and $\epsilon > 0$. Choose $n \in
{\mathbb N}$ very large and consider the map $\rho : E \to E \otimes
M_n({\mathbb C})$ given by
$$x \mapsto 1_E \otimes diag(0_{n-k}, \pi(x)) + x \otimes
diag(1_{n-k}, 0_k),$$ where $diag(0_{n-k}, \pi(x))$ is the block
diagonal element in $M_n({\mathbb C})$ whose first $n - k$ entries
down the diagonal are zero and the bottom block is given by $\pi(x)$,
while $diag(1_{n-k}, 0_k) \in M_n({\mathbb C})$ has $n - k$ $1$'s down
the diagonal followed by $k$ zeros. The key remarks about this choice
of connecting map are:

\begin{enumerate}
\item If $\frac{n - k}{n} > 1 - \epsilon$ then $|\tau \otimes tr_n
(\rho(x)) - \tau(x) | < 2\epsilon \|x\|$ for all $x \in E$.  That is,
in trace the connecting map $\rho$ is almost the same as the map $x
\mapsto x \otimes 1_{M_n}$ (which would be the natural connecting maps
to use if we were trying to construct $E \otimes \mathcal{U}$).

\item If $I \subset E$ is an ideal and there exists an element $x \in
I$ such that $\pi(x) \neq 0$ then the ideal generated by $\rho(I)$ is
all of $E\otimes M_n({\mathbb C})$.  This follows from the definition
of $\rho$ and the simplicity of $M_n({\mathbb C})$.  It is this fact
that will allow us to deduce simplicity in our inductive limits.

\item There exists a finite dimensional C$^*$-algebra $B \subset
E\otimes M_n({\mathbb C})$ with unit $e$ such that $e\rho(x) -
\rho(x)e = 0$ and $e\rho(x)e \in B$, for all $x \in E$.  (Let $B =
diag(0,\ldots,0,\pi(E)) \subset M_n({\mathbb C})$.) This remark will
immediately imply that our inductive limits satisfy the finite
dimensional approximation property which defines Popa algebras.

\item The representation $\pi \otimes \id : E \otimes M_n({\mathbb C})
\to M_k \otimes M_n({\mathbb C})$ is again a finite dimensional
representation and hence this whole procedure can be reapplied to the
algebra $E\otimes M_n({\mathbb C})$ (thus yielding an inductive
system).
\end{enumerate}

We now enter the gory details.  So let $E$ be a residually finite
dimensional C$^*$-algebra and $\pi_i : E \to M_{k(i)} ({\mathbb C})$
be a separating sequence of representations.  In fact, we will assume
that for every $x \in E$, $\| x \| = \lim_i \| \pi_i (x) \|$ (taking
direct sums, it is not hard to see that every residually finite
dimensional C$^*$-algebra has such a sequence).  Note that for every
$n \in {\mathbb N}$, $\pi_i \otimes id : E \otimes M_n({\mathbb C})
\to M_{k(i)} \otimes M_n({\mathbb C})$ is a separating sequence of the
same type.

Now choose natural numbers $1 = n(0) \leq n(1) \leq n(2) \leq \ldots$
such that $$\frac{n(0)n(1)\cdots n(j-1)k(j)}{n(j)} < 2^{-j},$$ for all
$j \in {\mathbb N}$.  One then defines algebras $E = E_0, E_1 = E_0
\otimes M_{n(1)}, E_2 = E_1 \otimes M_{n(2)}, E_3 = E_2 \otimes
M_{n(3)}, \ldots $ and inclusions $\rho_i : E_i \hookrightarrow E_{i +
1}$ as in the basic construction described above where the inclusion
$\rho_i$ uses the finite dimensional representation $\pi_{i+1} \otimes id
\otimes \cdots \otimes id : E \otimes M_{n(1)} \otimes \cdots \otimes
M_{n(i)} \to M_{k(i+1)} \otimes M_{n(1)} \otimes \cdots \otimes
M_{n(i)} $ in the lower right hand corner.  Letting $\Phi_{j,i} : E_i
\to E_j$, $i \leq j$, be defined by $\Phi_{j,i} = \rho_{j-1}\circ
\cdots \circ \rho_i$ we get an inductive system $\{ E_i, \Phi_{j,i}
\}$.

The are some projections in the above inductive system which we will
need.  Let $P_i \in E_i$ be the projection $$P_i = 1_{E_{i-1}} \otimes
diag(1_{n(i) - n(1)\cdots n(i-1)k(i)}, 0_{n(1)\cdots n(i-1)k(i)}).$$
Note that $P_{i+1}$ commutes with all of $\rho_i (E_i)$ (and, in
particular, with $\rho_i (P_i)$). Note also that if we write $E_{i+1}
= E_i \otimes M_{n(i+1)}$ then $$P_{i+1}\rho_i (P_i) = P_i \otimes
diag(1_{n(i+1) - n(1)\cdots n(i)k(i+1)}, 0_{n(1)\cdots n(i)k(i+1)}).$$

Letting $A$ be the inductive limit of the inductive system above, we
only have to show that $A$ is the Popa algebra we are after.

{\noindent\it Proof of Theorem \ref{thm:basicconstruction}:} We keep
all the notation above.  We leave it to the reader to verify that $A$
is a Popa algebra as this follows from our remarks above and the
construction of $A$. (That $A$ is unital and satisfies the right
finite dimensional approximation property is obvious while simplicity
follows from the remark that any ideal in A must eventually intersect
some $E_i$ (cf.\ \cite[Lemma III.4.1]{davidson}).)  Note also that $A$ was
constructed as an inductive limit of matrices over $E$ and hence
belongs to the class $\mathfrak{C}$ when $E$ does.

Now observe that given a trace $\tau \in T(E)$ we can define traces
$\tau_j \in T(E_j)$ by $\tau_j = \tau\otimes tr_{n(1)} \otimes \cdots
\otimes tr_{n(j)}$.  Then the embedding $\rho_j : E_j \to E_{j+1}$
almost intertwines $\tau_j$ and $\tau_{j+1}$.  More precisely, a
straightforward (but rather unpleasant) calculation shows that for $i
< j$, $$\tau_j (\Phi_{j,i} (x)) = \frac{\prod\limits_{s = i}^{j-1}
(n(s+1) - n(1)n(2)\cdots n(s)k(s+1))}{\prod\limits_{s = i}^{j-1}
n(s+1)} \tau_i (x) + \lambda_{i,j} \eta_{i,j}(x),$$ where
$\lambda_{i,j} = 1 - \frac{\prod\limits_{s = i}^{j-1} (n(s+1) -
n(1)n(2)\cdots n(s)k(s+1))}{\prod\limits_{s = i}^{j-1} n(s+1)}$ and
$\eta_{i,j}$ is some tracial state on $E_i$. Hence we get the estimate
$$|\tau_j (\Phi_{j,i}(x)) - \tau_i (x) | \leq 2\lambda_{i,j} \| x \|,$$
for all $x \in E_i$.  But, it can be shown by induction that
$\prod\limits_{s=i}^{j-1} (1 - \frac{n(1)n(2)\cdots
n(s)k(s+1)}{n(s+1)}) \geq \prod\limits_{s=i}^{j-1} (1 - 2^{-s-1}) \geq
1 - 2^{-i} + 2^{-j} \geq 1 - 2^{-i},$ for all $i < j \in {\mathbb
N}$.  Hence we get that $$|\lambda_{i,j}| = | 1 -
\prod\limits_{s=i}^{j-1} (1 - \frac{n(1)n(2)\cdots
n(s)k(s+1)}{n(s+1)}) | \leq 2^{-i}.$$

We have almost established part (1) in Theorem
\ref{thm:basicconstruction}. For each $i \in {\mathbb N}$, extend
$\tau_i$ to a state on $A$ (after identifying $E_i$ with it's image in
$A$).  It is clear that if we take any weak-$*$ cluster point,
$\gamma$, of this sequence then we will get a trace on $A$.  Moreover,
by the estimates above, we have that for each $x \in E_i$, $$|\gamma(x)
- \tau_i (x)| \leq 2^{-i} \| x \|.$$ Since we always have
$\tau$-preserving embeddings of $E$ into $E_i$, it should be clear how
to construct the embedding $\rho$ in the statement of the theorem.

Our last task is to prove that $\pi_{\gamma} (A)^{\prime\prime} \cong 
\pi_{\tau} (E)^{\prime\prime} \bar{\otimes} R$.  To do this, we will need 
to study the projections $P_j \in E_j$ defined above.  The idea is that 
we will use the $P_j$'s to construct different projections $Q^{(i)} \in 
\pi_{\gamma} (A)^{\prime\prime}$ with the following properties: 

\begin{enumerate}
\item $Q^{(i)} = P_i Q^{(i+1)} = Q^{(i+1)} P_i$ and hence $Q^{(i)}
\leq Q^{(i+1)}$ for all $i \in {\mathbb N}$.

\item $Q^{(i+1)} \in \pi_{\gamma} (E_i)^{\prime}$, for all $i \in
{\mathbb N}$ (where we have identified $E_i$ with it's image in $A$).

\item $\gamma(Q^{(i)}) \geq 1 - 2^{-i}$.

\item For each $i \in {\mathbb N}$,
$\frac{\gamma(Q^{(i+1)}x)}{\gamma(Q^{(i+1)})} = \tau_i (x)$, for all
$x \in E_i$.

\item The natural inclusion of the weak closure of $Q^{(i)}
\pi_{\gamma} (E_{i-1}) Q^{(i)}$ into the weak closure of
$Q^{(i+1)}\pi_{\gamma} (E_{i})Q^{(i+1)}$ (which is a natural inclusion
by (1) above) is isomorphic to the (non-unital) inclusion
$\pi_{\tau_{i-1}} (E_{i-1})^{\prime\prime} \hookrightarrow
\pi_{\tau_{i}}(E_{i})^{\prime\prime} \cong \pi_{\tau_{i-1}}
(E_{i-1})^{\prime\prime} \otimes M_{n(i)}$ given by $$x \mapsto
x\otimes diag(1_{n(i) - n(1)n(2)\cdots n(i-1)k(i)}, 0_{n(1)n(2)\cdots
n(i-1)k(i)}).$$
\end{enumerate}
 
We claim that the construction of such $Q^{(i)}$'s will complete the
proof.  Indeed, if we can do this then one uses part (5) and Elliott's
approximate intertwining argument to compare the (non-unital)
inclusions $Q^{(1)} \pi_{\gamma} (E_0) \subset Q^{(2)} \pi_{\gamma}
(E_1) \subset \ldots$ to the natural (unital) inclusions $E_0 \subset
E_0 \otimes M_{n(1)} \subset E_0 \otimes M_{n(1)} \otimes M_{n(2)}
\subset \ldots$.  Part (3) ensures that the former sequence recaptures
$\pi_{\gamma} (A)^{\prime\prime}$ while the latter sequence gives $E_0
\otimes \mathcal{U}$ in the limit, where $\mathcal{U}$ is a UHF
algebra, and hence the weak closure will be as desired and the proof
will be complete.

The construction of the $Q^{(i)}$'s is fairly simple.  For each $i \in
{\mathbb N}$ we define projections $Q^{(i)}_n = \pi_{\gamma}(P_i
P_{i+1} \cdots P_{i+n}) \in \pi_{\gamma}(A)^{\prime\prime}$.  Since
the $Q^{(i)}_n$'s are decreasing (as $n \to \infty$), there exists a
strong operator topology limit.  Define $Q^{(i)} = {\rm
sot}-\lim_{n\to \infty} Q^{(i)}_n$.  Then $Q^{(i)}$ is a projection
and it is straightforward to verify conditions (1) and (2)
above. (Recall that $P_j$ commutes with $\Phi_{j,i} (E_i)$ whenever $i
< j$.)  Thus we are left to verify the last three conditions.

Proof of (3). It suffices to show that $\gamma (Q^{(i)}_n) \geq 1 -
   2^{i} - 2^{-i-n}$ for all n.  But we may identify $Q^{(i)}_n$ with
   a projection in $E_{i+n}$ and so using the first part of the proof
   of this theorem (and using the identification) we get $|
   \gamma(Q^{(i)}_n) - \tau_{i+n}(Q^{(i)}_n)| < 2^{-i-n}$.  However,
   it follows from the construction of $Q^{(i)}_n$ that $$\tau_{i+n}
   (Q^{(i)}_n) = \prod\limits_{s = 1}^{n} tr_{n(i+s)}
   (diag(1_{n(i+s) - n(1)\cdots n(i+s-1)k(i+s)}, 0_{ n(1)\cdots
   n(i+s-1)k(i+s)})).$$ Thus, by the calculations given in the first
   part of the proof of this theorem, we see that $\tau_{i+n}
   (Q^{(i)}_n) \geq 1 - 2^{-i}$ and hence $\gamma(Q^{(i)}_n) \geq 1 -
   2^{i} - 2^{-i-n}$.

Proof of (4). We must show that if $x \in E_{i-1}$ then $\tau_{i-1}
    (x)\gamma(Q^{(i)}) = \lim_{n\to \infty} \gamma(Q^{(i)}_n x)$.  In
    order to show this, it suffices to prove that $$\tau_{i-1} (x) =
    \frac{\tau_{i+n}(Q^{(i)}_n x)}{\tau_{i+n}(Q^{(i)}_n)},$$ for all
    $n \in {\mathbb N}$.  This last equality, however, is evident from
    the construction.

Proof of (5). It suffices to show (essentially due to the uniqueness,
up to unitary equivalence, of GNS representations together with part
(4)) that there exists a surjective $*$-homomorphism $\eta : E_{i-1}
\otimes M_{n(i)} \to Q^{(i+1)}\pi_{\gamma} (E_i)$ such that for every
$x \in E_{i-1}$, $$\eta(x\otimes diag(1_{n(i) - n(1)\cdots
n(i-1)k(i)}, 0_{n(1)\cdots n(i-1)k(i)})) = Q^{(i)}\pi_{\gamma}(x) =
Q^{(i+1)}P_i \pi_{\gamma}(x).$$ But since
$$P_i \rho_{i-1}(x) = x\otimes diag(1_{n(i) - n(1)\cdots n(i-1)k(i)},
0_{n(1)\cdots n(i-1)k(i)}) \in E_i = E_{i-1} \otimes M_{n(i)},$$ we
get the desired homomorphism by identifying $E_{i-1} \otimes M_{n(i)}$
with it's image in $A$, passing to the GNS construction $\pi_{\gamma}$
and then cutting down by $Q^{(i+1)}$. \ \ $\square$

As mentioned above, the construction used in this section is a simple
adaptation of the Bratteli systems introduced in
\cite{dadarlat:nonnuclearsubalgebras}.  Chris Phillips is currently
writing a manuscript \cite{phillips:nonclassification} which contains
a detailed exposition of (a generalization of) these inductive
systems.  He also proves a number of results concerning the UCT
\cite{rosenberg-schochet} which we will find use for later.  The two
UCT results which are needed in this work are stated below.  See
section 2 of \cite{phillips:nonclassification} for a detailed
explanation of what is meant by the UCT and the proofs of the
following results.

\begin{subprop}[Two out of Three Principle]
\label{thm:twooutofthreeprinciple}
If $0 \to I \to E \to B \to 0$ is a semi-split extension 
(i.e.\ there exists a contractive c.p. splitting $B \to E$) 
and any two of $I$, $E$ and $B$ satisfy the UCT then so does 
the third.
\end{subprop}

\begin{subprop}
\label{thm:UCTinductivelimits}
Let $\phi_n : A_n \to A_{n+1}$ be injective $*$-homomorphisms and 
let $A$ denote the inductive limit of this sequence.  Assume further 
that for each $n$, there exists a contractive c.p. map $\psi_n : 
A_{n+1} \to A_n$ such that $\psi_n \circ \phi_n = id_{A_n}$ for 
all $n$.  If each $A_n$ satisfies the UCT then so does $A$.
\end{subprop}

\begin{subrem}
\label{thm:UCTremark}
Note that the inductive limit construction used in Theorem
\ref{thm:basicconstruction} does have the one-sided c.p. inverses
required in Proposition \ref{thm:UCTinductivelimits} and hence the
Popa algebra constructed will satisfy the UCT whenever the original
residually finite dimensional algebra does.  To see that such
c.p. inverses exist, we let $\pi: E \to M_k ({\mathbb C})$ be a
$*$-homomorphism and $\rho : E \to E\otimes M_n ({\mathbb C})$ be as
in the basic construction.  The desired map $E\otimes M_n ({\mathbb
C}) \to E \cong E \otimes e_{1,1}$ is just given by compressing to the
(1,1) corner. \end{subrem}

\subsection{Miscellaneous}
Here we collect a few facts for future reference.  The first two are
well known and the second two are just a bit of trickery.  We begin
with a simple adaptation of the fact that if $M$ is a von Neumann
algebra with faithful, normal tracial state $\tau$ and $1_M \in N
\subset M$ is a sub-von Neumann algebra then there always exists a
$\tau$-preserving (hence faithful) conditional expectation $M \to N$
(cf.\ \cite[Exercise 8.7.28]{kadison-ringrose}).

\begin{sublem}
\label{thm:tracepreservingconditionalexpectation}
Let $A$ be a C$^*$-algebra, $\tau \in \TA$ be a  tracial state 
and $1_A \in B \subset A$ be a finite dimensional subalgebra. Then 
there exists a conditional expectation $\Phi_B : A \to B$ such that 
$\tau \circ \Phi_B = \tau$.
\end{sublem}

\begin{proof} Assume first that $\tau|_B$ is faithful. 
Write $B \cong M_{n(1)} ({\mathbb C}) \oplus \cdots \oplus M_{n(k)}
({\mathbb C})$, and let $\{ e^{(1)}_{i,j} \}_{1 \leq i,j \leq n(1)} \cup
\ldots \cup \{ e^{(k)}_{i,j} \}_{1 \leq i,j \leq n(k)}$ be a system of
matrix units for $B$. Then the desired conditional expectation is
given by $$\Phi_B (x) = \sum\limits_{s = 1}^{k} \sum\limits_{i,j =
1}^{n(s)} \frac{\tau(xe_{i,j}^{(s)})}{\tau(e_{i,i}^{(s)})}
e_{j,i}^{(s)}.$$ 

When $\tau|_B$ is not faithful, the formula above no longer makes
sense.  However, one can decompose $B$ as the direct sum of two finite
dimensional algebras, $B_0 \oplus B_f$, where $\tau|_{B_0} = 0$ and
$\tau|_{B_f}$ is faithful.  Letting $e_0$ (resp.\ $e_f$) be the unit of
$B_0$ (resp.\ $B_f$) we get a $\tau$-preserving conditional
expectation by mapping each $a \in A$ to $E_{B_0} (e_0 a e_0) +
E_{B_f} (e_f a e_f)$, where $E_{B_0} : e_0 A e_0 \to B_0$ is any
conditional expectation (which exist by finite dimensionality) and
$E_{B_f} : e_f A e_f \to B_f$ is a $\tau|_{e_f A e_f}$ preserving
conditional expectation as in the first part of the proof.
\end{proof}

\begin{sublem}
\label{thm:tracepreservingembedding}
If $B$ is a finite dimensional C$^*$-algebra with tracial state $\tau$
then for every $\varepsilon > 0$ there exists  $n \in {\mathbb
N}$ and a unital $*$-monomorphism $\rho : B \hookrightarrow
M_n({\mathbb C})$ such that $| \tau (x) - tr_n \circ \rho (x) | <
\varepsilon \| x \|$ for every $x \in B$.
\end{sublem}

\begin{proof} 
If $\tau$ is a {\em rational} convex combination of extreme traces
then one can find an honestly trace preserving embedding by inflating
the summands of $B$ according to the rational numbers appearing in the
convex combination. The general case then follows by approximation.
\end{proof}

We remind the reader of the following theorem of Voiculescu (cf.\
\cite[Theorem 1]{dvv:QDhomotopy}): A separable, unital C$^*$-algebra
$A$ is quasidiagonal (QD) if and only if there exists a sequence of
u.c.p. maps $\varphi_n : A \to M_{k(n)} ({\mathbb C})$ which are
asymptotically multiplicative (i.e.\ $\| \varphi_n (ab) -
\varphi_n(a)\varphi_n(b) \| \to 0$ for all $a,b \in A$) and
asymptotically isometric (i.e.\ $\|a \| = \lim \| \varphi_n (a) \|$,
for all $a \in A$).

\begin{sublem}
\label{thm:trickery}
Assume that $A$ is a QD C$^*$-algebra and $\psi_n : A \to M_{l(n)}
({\mathbb C})$ is an asymptotically multiplicative (but not
necessarily asymptotically isometric) sequence of u.c.p. maps. Then
there exists a sequence of u.c.p. maps $\Phi_n : A \to M_{t(n)}
({\mathbb C})$ which are asymptotically multiplicative, asymptotically
isometric and such that $| tr_{t(n)}(\Phi_n (a)) - tr_{l(n)} (\psi_n
(a)) | \to 0$ for all $a \in A$.
\end{sublem}

\begin{proof}
Let $\varphi_n : A \to M_{k(n)} ({\mathbb C})$ be an asymptotically
multiplicative, asymptotically isometric sequence of u.c.p. maps.
Choose integers $s(n)$ such that $\frac{s(n)l(n)}{s(n)l(n) + k(n)} \to
1$ (i.e.\ such that $\frac{k(n)}{s(n)} \to 0$).  Then one defines
$\Phi_n : A \to M_{s(n)l(n) + k(n)} ({\mathbb C})$ to be the block
diagonal map with one summand equal to $\varphi_n$ and $s(n)$ summands
equal to $\psi_n$.
\end{proof}

Note that the previous lemma  can also be formulated in terms of
$*$-homomorphisms when $A$ is a residually finite dimensional
$C^*$-algebra.

Though we will try to keep everything unital, non-unital maps are sometimes 
unavoidable.  The next lemma keeps everything running smoothly. 

\begin{sublem}
\label{thm:nonunital}
Let $\phi_n : A \to M_{k(n)} ({\mathbb C})$ be contractive c.p.  maps
(though we still assume $A$ is unital) which are asymptotically
multiplicative with respect to the 2-norms on $M_{k(n)} ({\mathbb
C})$.  If $tr_{k(n)} (\phi_n (1_A)) \to 1$ then there exist
u.c.p. maps $\psi_n : A \to M_{k(n)} ({\mathbb C})$ which are also
asymptotically multiplicative (with respect to 2-norms) and such that
$| tr_{k(n)}(\psi_n (a)) - tr_{k(n)}(\phi_n(a)) | \to 0$ as $n \to
\infty$, with convergence being uniform on the unit ball of $A$.
\end{sublem}

\begin{proof} By \cite[Lemma 2.2]{choi-effros:injectivityandoperatorspaces}
we can find u.c.p. maps $\psi_n : A \to M_{k(n)} ({\mathbb C})$ such
that $\phi_n (a) = c_n \psi_n (a) c_n$ for all $n$ and $a \in A$,
where $c_n = \phi_n (1_A)^{1/2}$.  Since $tr_{k(n)} (\phi_n (1_A)) \to
1$ it follows that $\| c_n - 1_{M_{k(n)}} \|_2 \to 0$.  Using this
remark, the Cauchy-Schwartz inequality and the general inequality $\|
x - c_n x c_n \|_2 \leq 2 \| x \| \| c_n - 1_{M_{k(n)}} \|_2$ (for all
$x \in M_{k(n)} ({\mathbb C})$) it is straightforward to verify that
$\psi_n$ have the properties asserted in the statement of the lemma.
\end{proof}

\section{Approximately Finite Dimensional Traces}

We remind the reader that all C$^*$-algebras (resp.\ von Neumann
algebras) in these notes are assumed to be separable and unital
(resp.\ have separable preduals).  We also remind the reader that all
notation not defined below should have been explained in Section 2.1.

We now define the tracial invariants mentioned in the introduction. 
We will be concerned with those traces that can be approximated by
almost multiplicative maps to matrix algebras.  We can require these
maps to be almost multiplicative either in the operator norm or the
2-norm and we can require either weak-$*$ or norm convergence in the
dual space.  Hence we get four types of approximation.  Kirchberg has
studied similar sets of traces in his work on Connes' embedding
problem (see \cite[Lemma 4.5]{kirchberg:invent}) and residually finite
groups with Kazhdan's property T (see \cite[Definition
3.1]{kirchberg:propertyTgroups}).  Indeed, it is not hard to see that
our definition of $\TAwafd$ below agrees with Kirchberg's notion of
{\em liftable} traces from \cite{kirchberg:propertyTgroups}.  The main
theorem concerning the set $\TAwafd$ (Theorem
\ref{thm:mainthmsection3} below) will thus depend heavily on
\cite[Proposition 3.2]{kirchberg:propertyTgroups}.

\begin{defn} 
\label{thm:AFDtraces}
Let $A$ be a C$^*$-algebra and $\TA$ denote the (possibly empty) set
of tracial states on $A$.  We will say that a trace $\tau \in \TA$ is
{\em approximately finite dimensional} if there exists a sequence of
u.c.p. maps $\phi_n : A \to M_{k(n)} (\mathbb C)$ such that $\| \phi_n
(ab) - \phi_n (a) \phi_n (b) \| \to 0$ and $\tau (a) = \lim_{n \to
\infty} tr_{k(n)} \circ \phi_n (a)$, for all $a,b \in A$.  If we have
norm convergence in the dual space $A^*$ (i.e.\ $\| \tau - tr_{k(n)}
\circ \phi_n \|_{A^*} \to 0$) then $\tau$ will be called {\em
uniformly approximately finite dimensional}.

A trace $\tau \in \TA$ is called {\em weakly approximately finite
dimensional} if there exists a sequence of u.c.p. maps $\phi_n : A \to
M_{k(n)} (\mathbb C)$ such that $\| \phi_n (ab) - \phi_n (a) \phi_n
(b) \|_2 \to 0$ and $\tau (a) = \lim_{n \to \infty} tr_{k(n)} \circ
\phi_n (a)$, for all $a,b \in A$.  Similarly, $\tau$ will be called
{\em uniformly weakly approximately finite dimensional} if
$tr_{k(n)}\circ\phi_n \to \tau$ in the norm topology on $A^*$.
\end{defn}

We then put 

$$\TAafd = \{ \tau \in \TA : \tau {\rm \ is \ approximately \ finite \
dimensional} \},$$

$$\TAwafd = \{ \tau \in \TA : \tau {\rm \ is \
weakly \ approximately \ finite \ dimensional} \}, $$

$$\UTAafd = \{ \tau \in \TA : \tau {\rm \ is \ uniformly \ approximately
\ finite \ dimensional} \}, $$

$$\UTAwafd = \{ \tau \in \TA : \tau {\rm \ is \ uniformly \ 
weakly \ approximately \ finite \ dimensional} \}.$$ 

Evidently we have the following inclusions of the sets defined above. 

$$\begin{array}{ccccc}
 \TA & \supset  & \TAwafd  & \supset & \TAafd \\
     &          &    \cup  &         & \cup \\
     &          & \UTAwafd & \supset & \UTAafd 
\end{array}$$

We will see in Section 6 that if $A$ is locally reflexive (e.g.\
nuclear or exact) then $\TA \supset \TAwafd = \UTAwafd \supset \TAafd
= \UTAafd$.

Though they may seem unnatural at first, the definitions above are
inspired by the theories of quasidiagonal and nuclear C$^*$-algebras.

Recall that a C$^*$-algebra $A$ is quasidiagonal (QD) if
there exists a faithful representation $\pi : A \to B(H)$
 such that one can find an increasing sequence
of finite rank projections $P_1 \leq P_2 \leq \ldots$ with the
property that $\| \pi(a)P_n - P_n \pi(a) \| \to 0$ for all $a \in A$
and $P_n \to 1_H$ in the strong operator topology. (This is not the 
right definition for non-separable algebras.)

\begin{example}(cf.\ \cite[2.4]{dvv:QDsurvey}, 
\cite[Proposition 6.1]{brown:QDsurvey}) If $A$ is  QD then $\TAafd
\neq \emptyset$ (in the non-unital case it can happen that $\TA =
\emptyset$; e.g.\ the suspension of a Cuntz algebra).  Indeed if $\pi
: A \to B(H)$ and $P_1 \leq P_2 \leq \ldots$ are as in the definition
of quasidiagonality then one defines u.c.p. maps by $$\phi_n(a) = P_n
\pi(a) P_n.$$ The asymptotic commutativity of $P_n$ ensures that
$\phi_n$ are asymptotically multiplicative in norm.  Note also that $P_n B(H)
P_n \cong M_{rank(P_n)} (\mathbb C)$.  Finally, a straightforward
calculation shows that any weak-$*$ cluster point of the sequence of
states $\{ tr_{rank(P_n)} \circ \phi_n \}$ is necessarily a tracial
state and hence $\TAafd \neq \emptyset$.
\end{example}

From the previous example we see that the approximately finitely
dimensional traces form a very natural subset of the tracial space of
a QD C$^*$-algebra.  Will we see later that $\UTAafd$ is also fairly
natural (at least in the classification program).  A much deeper fact
is that the sets $\TAwafd$ and $\UTAwafd$ are large for an important
class of C$^*$-algebras (see also Corollary \ref{thm:WEPcase}).

\begin{prop}
\label{thm:nuclearcase}
If $A$ is a nuclear C$^*$-algebra then $\TA = \TAwafd = \UTAwafd$.
\end{prop}

\begin{proof}
The proof is a simple consequence of the following deep fact: If $A$
is nuclear and $\tau \in \TA$ then $\pi_{\tau}(A)^{\prime\prime}$ is
hyperfinite (cf.\ \cite{choi-effros:nuclearityandinjectivity},
\cite{connes:classification},
\cite{popa:injectiveimplieshyperfiniteI},
\cite{popa:injectiveimplieshyperfiniteII} ).  Given this result, the
equation above follows easily from Lemmas
\ref{thm:tracepreservingconditionalexpectation} and
\ref{thm:tracepreservingembedding}.
\end{proof}

\begin{prop}
The sets $\TAafd$, $\TAwafd$, $\UTAafd$ and $\UTAwafd$ are all convex.
The sets $\TAafd$ and $\TAwafd$ are closed in the weak-$*$ topology
while $\UTAwafd$ and $\UTAafd$ are closed in norm (i.e.\ the norm on
$A^*$) and thus, by the Hahn-Banach theorem, closed in the weak
topology coming from $A^{**}$.
\end{prop}

\begin{proof} It is not hard to verify that the first assertion 
follows from Lemma \ref{thm:tracepreservingembedding}.  Since $\TAafd$
and $\TAwafd$ (resp.\ $\UTAwafd$ and $\UTAafd$) are defined via
weak-$*$ convergence (resp.\ norm convergence), it is also easy to see
that these sets are closed in this topology.
\end{proof}

We will soon see that if $\Gamma$ is a non-amenable, residually
finite, discrete group then the canonical trace on $C^* (\Gamma)$
(which gives the left regular representation in the GNS construction)
is always a weak-$*$ limit of uniformly approximately finite
dimensional traces, but is not itself uniformly weakly approximately
finite dimensional.  Thus, the sets $\UTAwafd$ and $\UTAafd$ need not
be weak-$*$ closed in general.

In \cite[Definition 3.1]{kirchberg:propertyTgroups} Kirchberg
introduced the notion of a {\em liftable} trace.  His definition is as
follows: $\tau \in \TA$ is {\em liftable} if there exists a completely
positively liftable $*$-homomorphism $\pi : A \to R^{\omega}$ such
that $\tau_{\omega}\circ\pi = \tau$. (Here, $R^{\omega}$ is as defined
in Section 2.1 and `completely positively liftable' means that there
exists a completely positive lifting (from $A \to l^{\infty}(R)$) for
$\pi$.)  To the experts it is immediate that this notion is the same
as weakly approximately finite dimensional, as defined above, and it
is a good exercise for those not familiar with these notions.  In any
case, the next proposition is now a direct consequence of \cite[Lemma
3.4]{kirchberg:propertyTgroups}.

\begin{prop} 
$\TAwafd$ is a face in $\TA$.
\end{prop}

We will soon see that $\UTAwafd$ is also a face, but we do not
 know whether the sets $\TAafd$ and $\UTAafd$ are always faces.

The hardest part in the proof of the following theorem is the
implication $(7) \Longrightarrow (8)$.  As previously indicated, this
was already done by Kirchberg in \cite[Proposition
3.2]{kirchberg:propertyTgroups}.  His proof amounts to a very clever
(and technically difficult) reduction to \cite[Theorem
1.2.2]{connes:classification}.  Indeed, as will become clear, most of
the content of this theorem can be traced back to the fundamental work
of Alain Connes \cite{connes:classification}.

\begin{thm}
\label{thm:mainthmsection3}

Let $A \subset B(H)$ be  in general
position (i.e.\ $A \cap \mathcal{K}(H) = \{0 \}$) and $\tau \in \TA$.  
Then the following are equivalent:

\begin{enumerate}
\item $\tau \in \TAwafd$.

\item There exist finite rank projections $P_n \in B(H)$ (not
necessarily increasing) such that $$\frac{\| aP_n - P_n a \|_{HS}}{\|
P_n \|_{HS}} \to 0 \ {\rm and \ } \tau(a) = \lim_{n \to \infty}
\frac{<aP_n, P_n>_{HS}}{<P_n, P_n>_{HS}},$$ for all $a \in A.$

\item For all $a_1, \ldots, a_n, b_1, \ldots, b_n \in A$, $$|
\tau(\sum\limits_{i = 1}^{n} a_i b_i^*) | \leq \| \sum\limits_{i =
1}^{n} a_i \otimes b_i^* \|_{A\otimes A^{op}}.$$

\item For all $a_1, \ldots, a_n, b_1, \ldots, b_n \in A$, $$\|
\pi_{\tau} \odot \pi_{\tau}^{op} (\sum\limits_{i = 1}^{n} a_i \otimes
b_i) \| \leq \| \sum\limits_{i = 1}^{n} a_i \otimes b_i \|_{A\otimes
A^{op}}. $$

\item $\pi_{\tau} \odot \pi_{\tau}^{op}$ extends to a representation
of $A \otimes A^{op}$.

\item There exists a u.c.p. map $\Phi : B(H) \to \pi_{\tau} (A)^{\prime
\prime}$ such that $\Phi(a) = \pi_{\tau} (a)$, for all $a \in A$.

\item $\tau$ extends to a hypertrace (i.e.\ there exists a state
$\varphi$ on $B(H)$ such that $\varphi|_A = \tau$ and $A \subset
B(H)_{\varphi} = \{ T \in B(H) : \varphi(ST) = \varphi(TS), \ {\rm for
\ all \ } S \in B(H)\}$).

\item There exists a $*$-monomorphism $\Phi : \pi_{\tau}
(A)^{\prime\prime} \hookrightarrow R^{\omega}$ such that
$\tau^{\prime\prime} = \tau_{\omega}\circ\Phi$ and $\Phi \circ
\pi_{\tau} : A \to R^{\omega}$ can be lifted to a u.c.p. map $A \to
l^{\infty}(R)$, where $\tau^{\prime\prime}$ is the vector trace on
$\pi_{\tau} (A)^{\prime\prime}$ induced by $\tau$.
\end{enumerate}
\end{thm}

\begin{proof}
(1) $\Longrightarrow$ (2). This follows from Voiculescu's
Theorem.  

(2) $\Longrightarrow$ (3) and (3) $\Longrightarrow$ (4) are
essentially the same as (6) $\Longrightarrow$ (5) and (5)
$\Longrightarrow$ (4), respectively, from \cite[Theorem
5.1]{connes:classification}.

(4) $\Longrightarrow$ (5) is immediate.

(5) $\Longrightarrow$ (6) is a consequence of \cite[Observation
3.0]{kirchberg:invent} but we remind the reader of Kirchberg's elegant
proof.  Since $A\otimes A^{op} \subset B(H) \otimes A^{op}$ we can
extend $\pi_{\tau} \otimes \pi_{\tau}^{op} : A\otimes A^{op} \to
B(L^2(A,\tau))$ to a completely positive map $\Phi : B(H) \otimes
A^{op} \to B(L^2(A,\tau))$.  Since $\Phi|_{A\otimes A^{op}}$ is a
homomorphism it follows that $A \otimes A^{op}$ (and, in particular,
$1\otimes A^{op}$) is in the multiplicative domain of $\Phi$.  Hence,
for every $T \in B(H)$, it follows that $\Phi(T \otimes 1) \in
\Phi(1\otimes A^{op})^{\prime} = \pi_{\tau}^{op} (A^{op})^{\prime} =
\pi_{\tau}(A)^{\prime\prime}$.

(6) $\Longrightarrow$ (7). Since $A$ is contained in the
multiplicative domain of $\Phi$, it is easy to verify that $\varphi(T)
= <\Phi(T)\eta_{\tau}, \eta_{\tau}>$ gives the desired hypertrace
(here, $\eta_{\tau} = \hat{1_A} \in L^2(A,\tau)$ is the canonical
trace vector).

(7) $\Longrightarrow$ (8). This follows from $(iii)$ $\Longrightarrow$
$(ii)$ in \cite[Proposition 3.2]{kirchberg:propertyTgroups}.  Indeed
Kirchberg's result states that, assuming (7), we can find a
u.c.p. liftable $*$-homomorphism $\Psi : A \to R^{\omega}$ such that
$\tau = \tau_{\omega} \circ \Psi$.  Hence the weak closure of
$\Psi(A)$ (inside $R^{\omega}$) will be canonically isomorphic to
$\pi_{\tau} (A)^{\prime\prime}$ and the desired map $\Phi$ is just
this identification.  

Finally, (8) $\Longrightarrow$ (1) is immediate to the reader who has
verified that liftable traces and weakly approximately finite
dimensional traces are the same thing.
\end{proof}

\begin{rem}
Note that the proofs of (5) $\Longrightarrow$ (6) and (6)
$\Longrightarrow$ (7) never used our assumption that $A$ is in general
position.  Hence {\em every weakly approximately finite dimensional
trace extends to a hypertrace in any faithful representation of $A$}.
\end{rem}

The space $\UTAwafd$ also admits a number of nice characterizations.
We thank Yasuyuki Kawahigashi, Sergei Neshveyev and Narutaka Ozawa for
discussions in Oberwolfach which slightly shortened our original proof
of (1) $\Longrightarrow$ (2).  Ozawa also added (6) to the list below
and showed us the elegant proof. 

\begin{thm}
\label{thm:mainthmUTAwafd}

For a trace $\tau \in \TA$, the following are equivalent:

\begin{enumerate}
\item $\tau \in \UTAwafd$.

\item There exist u.c.p. maps $\psi_n : A^{**} \to M_{k(n)} ({\mathbb
C})$ such that for each free ultrafilter $\omega \in \beta{\mathbb N}
\backslash {\mathbb N}$ we have $$\lim\limits_{n \to \omega} \| \psi_n
(xy) - \psi_n(x) \psi_n (y) \|_{2} = 0 \ {\rm and} \ \lim\limits_{n
\to \omega} tr_{k(n)} \circ \psi_n (x) = \tau^{**} (x),$$ for all $x,y
\in A^{**}$, where $\tau^{**}$ is the normal trace on $A^{**}$ induced
by $\tau$.

\item There exist u.c.p. maps $\psi_n : \pi_{\tau} (A)^{\prime\prime}
\to M_{k(n)} ({\mathbb C})$ such that for each free ultrafilter
$\omega \in \beta{\mathbb N} \backslash {\mathbb N}$ we have
$$\lim\limits_{n \to \omega} \| \psi_n (xy) - \psi_n(x) \psi_n (y)
\|_{2} = 0 \ {\rm and} \ \lim\limits_{n \to \omega} tr_{k(n)} \circ
\psi_n (x) = \tau^{\prime\prime} (x),$$ for all $x,y \in \pi_{\tau}
(A)^{\prime\prime}$, where $\tau^{\prime\prime}$ is the normal trace
on $\pi_{\tau} (A)^{\prime\prime} $ induced by $\tau$.

\item There exists a u.c.p. liftable, normal $*$-monomorphism $\sigma :
\pi_{\tau} (A)^{\prime\prime} \hookrightarrow R^{\omega}$.

\item $\pi_{\tau} (A)^{\prime\prime}$ is hyperfinite.

\item $\pi_{\tau} : A \to \pi_{\tau} (A)^{\prime\prime}$ is weakly
nuclear (i.e.\ there exists u.c.p.\ maps $\phi_n : A \to M_{k(n)}
({\mathbb C})$, $\psi_n : M_{k(n)} ({\mathbb C}) \to \pi_{\tau}
(A)^{\prime\prime}$ such that $\psi_n \circ \phi_n (a) \to \pi_{\tau}
(a)$ in the $\sigma$-weak topology for every $a \in A$).
\end{enumerate}
\end{thm}

\begin{proof} (1) $\Longrightarrow$ (2). Let $\phi_n : A \to M_{k(n)} 
({\mathbb C})$ be 2-norm asymptotically multiplicative u.c.p. maps
such that $tr_{k(n)} \circ \phi_n \to \tau$ in the norm on $A^*$.  Let
$\phi_n^{**} : A^{**} \to M_{k(n)} ({\mathbb C})$ be the canonical
extensions to the double dual.  Since $\| tr_{k(n)} \circ \phi_n -
\tau \|_{A^*} \to 0$, it follows that $tr_{k(n)} \circ \phi_n^{**} (x)
\to \tau^{**}(x)$ for every $x \in A^{**}$.  Hence we can construct a
u.c.p. map $\Phi : A^{**} \to R^{\omega}$ such that $i)$
$\tau_{\omega} \circ \Phi = \tau^{**}$, $ii)$ $\Phi|_{A}$ is a
$*$-homomorphism and $iii)$ $\Phi$ has a u.c.p. lifting $A^{**} \to
l^{\infty} (R)$.  If we knew that $\Phi$ was a homomorphism then it
would follow that $\lim\limits_{n \to \omega} \| \phi_n^{**} (xy) -
\phi_n^{**} (x) \phi_n^{**} (y) \|_{2} = 0$ for all $x,y \in A^{**}$
and hence this is what we will show.

First note that $\Phi$ is normal: if $\{ x_{\lambda} \} \subset
A^{**}_{sa}$ is a norm bounded, increasing net of self adjoint
elements with strong operator topology limit $x$ then $\{
\Phi(x_{\lambda})\}$ is increasing up to $\Phi(x)$ (in the strong
operator topology -- i.e.\ 2-norm) since $\Phi(x_{\lambda}) \leq
\Phi(x)$ and $\tau_{\omega}\circ (\Phi(x_{\lambda})) = \tau^{**}
(x_{\lambda}) \to \tau^{**} (x) = \tau_{\omega} \circ \Phi(x)$.  It
follows that $\Phi$ is continuous from the $\sigma$-weak topology on
$A^{**}$ to the $\sigma$-weak topology on $R^{\omega}$ (i.e.\ w.r.t.\
the weak-$*$ topologies coming from the preduals).  Letting $\Psi :
A^{**} \to R^{\omega}$ be the (weak-$*$ continuous) $*$-homomorphism
which extends $\Phi|_{A}$ (and which exists by universality of
$A^{**}$) it follows that $\Phi = \Psi$ since they are continuous and
agree on $A$.  Hence $\Phi$ is also multiplicative.

(2) $\Longrightarrow$ (3).  Assuming (2), we can use the maps $\psi_n$
to construct a u.c.p. liftable, $*$-homomorphism $\sigma : A^{**} \to
R^{\omega}$ such that $\tau_{\omega} \circ \sigma = \tau^{**}$.  It
follows that $\sigma(A^{**}) \cong \pi_{\tau} (A)^{\prime\prime}$ and,
hence, $A^{**} \cong ker(\sigma) \oplus \pi_{\tau}
(A)^{\prime\prime}$.  Restricting the maps $\psi_n$ to this non-unital
copy of $\pi_{\tau} (A)^{\prime\prime}$ gives c.p. maps with the
desired properties.  Then applying Lemma \ref{thm:nonunital} we can
replace these nonunital maps with unital ones and we get (3).

(3) $\Longrightarrow$ (4)  is immediate.

(4) $\Longrightarrow$ (5).  Since $l^{\infty}(R)$ is injective and there 
exists a (trace preserving) conditional expectation $R^{\omega} \to 
\sigma( \pi_{\tau} (A)^{\prime\prime})$, it is not hard to see that $ 
\pi_{\tau} (A)^{\prime\prime}$ must be injective, assuming (4).  

(5) $\Longrightarrow$ (1) is contained in the  proof of Proposition 
\ref{thm:nuclearcase} and hence (1) - (5) are equivalent.

(5) $\Longrightarrow$ (6) is trivial and hence we are left to prove
(6) $\Longrightarrow$ (5).  So assume that $\pi_{\tau} : A \to
\pi_{\tau} (A)^{\prime\prime}$ is weakly nuclear and $\phi_n : A \to
M_{k(n)} ({\mathbb C})$, $\psi_n : M_{k(n)} ({\mathbb C}) \to
\pi_{\tau} (A)^{\prime\prime}$ are u.c.p.\ maps whose composition
converges to $\pi_{\tau}$ in the point-$\sigma$-weak topology. Using
these maps it is not hard to see that the canonical homomorphism $A
\odot \pi_{\tau} (A)^{\prime} \to B(L^2(A,\tau))$, $a\otimes x \mapsto
\pi_{\tau} (a)x$, is continuous with respect to the minimal tensor product
norm. (Use the fact that the natural map on the maximal tensor product
approximately factorizes through $M_{k(n)} \otimes_{max} \pi_{\tau}
(A)^{\prime} = M_{k(n)} \otimes \pi_{\tau} (A)^{\prime}$ and hence
factors through the minimal tensor product.)  As in the proof of (5)
$\Longrightarrow$ (6) from the last theorem, it follows that there
exists a conditional expectation $B(L^2(A,\tau)) \to \pi_{\tau}
(A)^{\prime}$ and hence $\pi_{\tau} (A)^{\prime}$ is injective.  This
implies that $\pi_{\tau} (A)^{\prime\prime}$ is also injective and the
proof is complete.
\end{proof}

Note that part (3) in the previous theorem could be used as an
 abstract (i.e.\ representation free) definition of quasidiagonality,
 analogous to Voiculescu's abstract characterization \cite[Theorem
 1]{dvv:QDhomotopy}, in the setting of tracial von Neumann algebras.
 The equivalence of (3) and (5) would then say that quasidiagonality
 is equivalent to hyperfiniteness.  (Compare with the C$^*$-case where
 Dadarlat has constructed non-nuclear Popa algebras
 \cite{dadarlat:nonnuclearsubalgebras}.)

\begin{cor}
$\UTAwafd$ is always a face in $\TA$.
\end{cor}

\begin{proof} If $0 < s < 1$, $\tau, \gamma \in T(A)$ and $s\tau + 
(1-s)\gamma \in \UTAwafd$ then we can find (non-unital) normal
embeddings of $\pi_{\tau}(A)^{\prime\prime}$ and
$\pi_{\gamma}(A)^{\prime\prime}$ into the weak closure of the GNS
representation of $A$ with respect to $s\tau + (1-s)\gamma$ (cf.\
\cite[Proposition 3.3.5]{pedersen:book}).  As the latter algebra is
injective, it follows that both $\pi_{\tau}(A)^{\prime\prime}$ and
$\pi_{\gamma}(A)^{\prime\prime}$ are injective (hence hyperfinite).
Thus, both $\tau$ and $\gamma$ belong to $\UTAwafd$.
\end{proof}

The set $\TAafd$ seems a bit hard to get a handle on.  This is largely 
due, at least in this author's opinion, to our present lack of 
understanding of the class of QD C$^*$-algebras.  However, 
if one knows a C$^*$-algebra to be QD then the set $\TAafd$ is 
precisely the set of traces which can be encoded in the definition 
of quasidiagonality. 

\begin{prop}
\label{thm:QDcase}
Let $A \subset B(H)$ be in general position.  If $A$ is QD then there 
exists an increasing sequence of finite rank projections $P_1 \leq 
P_2 \leq \ldots$, converging strongly to the identity, which 
asymptotically commutes (in norm) with every element in $A$ and such 
that for each $\tau \in \TAafd$ there exists a subsequence $\{ n_k \}$ 
such that $$\frac{<aP_{n_k}, P_{n_k}>_{HS}}{<P_{n_k}, P_{n_k}>_{HS}} 
\to \tau(a), \ as \ k \to \infty,$$ for all $a \in A$.
\end{prop}

\begin{proof}
We claim that it suffices to prove that for every finite 
set $\mathfrak{F} \subset A$, finite dimensional subspace $X \subset H$,  
$\varepsilon > 0$ and trace $\tau \in \TAafd$ there exists a finite rank 
projection $P \in B(H)$ such that 
\begin{enumerate}
\item $\| [a,P] \| < \varepsilon$ for all $a \in \mathfrak{F}$.

\item $P(x) = x$ for all $x \in X$. 

\item $| \frac{<aP, P>_{HS}}{<P, P>_{HS}} - \tau(a) | < \varepsilon$ 
for all $a \in \mathfrak{F}$.
\end{enumerate}

Assume for the moment that we were able to prove this local version.
Then, if $\{ a_n \} \subset A$ is a sequence which is dense in the
unit ball of $A$ and $\{ \tau_j \}$ is any sequence of traces in
$\TAafd$ we could apply the above local approximation property to
construct a sequence $P_1 \leq P_2 \leq \ldots$ which was converging
strongly to the identity, asymptotically commuting in norm with $A$
and such that $$| \frac{<a_i P_n, P_n>_{HS}}{<P_n, P_n>_{HS}} -
\tau_n(a_i) | < 1/n$$ for all $n \in {\mathbb N}$ and $1 \leq i \leq
n$.  Since $A$ is separable, the weak-$*$ topology on $\TA$ is
metrizable and hence we can always find a sequence of traces $\{
\tau_j \} \subset \TAafd$ such that there exists a weak-$*$ dense
subset $Y \subset \TAafd$ with the property that every element of $Y$
appears infinitely many times in the sequence $\{ \tau_j \}$.  The
sequence of projections associated with such a sequence of traces will
have all the properties asserted in the statement of the
proposition. Hence it suffices to prove the local statement in the
first paragraph of the proof.

The required local statement is now a consequence of Voiculescu's
Theorem (version \ref{thm:technicalVoiculescuThm}) and a little
trickery.  Let $\tau \in \TAafd$ be arbitrary. Since $A$ is QD, by
Lemma \ref{thm:trickery} we can find a sequence of u.c.p.  maps
$\phi_n : A \to M_{k(n)}({\mathbb C})$ which are asymptotically
multiplicative, asymptotically isometric and such that $tr_{k(n)}
\circ \phi_n \to \tau$ in the weak-$*$ topology.  If a finite set
$\mathfrak{F} \subset A$ and $\epsilon > 0$ are given then, by passing
to a subsequence if necessary, we may assume that $\| \phi_n (ab) -
\phi_n(a)\phi_n(b) \| < \varepsilon$ and $|tr_{k(n)} \circ \phi_n (a)
- \tau(a) | < \varepsilon$ for all $n$ and for all $a,b \in
\mathfrak{F}$.  Letting $K = \oplus_n {\mathbb C}^{k(n)}$ and $\Phi =
\oplus_n \phi_n : A \to B(K)$ we have that $\Phi$ is a faithful
$*$-homomorphism modulo the compacts.  Hence we can find a unitary
operator $U : K \to H$ such that $U \Phi(a) U^*$ is nearly equal (in
norm) to $a$, for all $a \in \mathfrak{F}$.  Hence, if $Q_s \in B(K)$
is the orthogonal projection onto $\oplus_1^s {\mathbb C}^{k(n)}$ we
have that $\| [UQ_s U^*, a]\|$ is small for all $s \in {\mathbb N}$
and for all $a \in A$.  Moreover, compressing $a \in \mathfrak{F}$ to
the range of $UQ_s U^*$ will almost recover the trace $\tau$ (for all
$s$).  Finally, if $X \subset H$ is any finite dimensional subspace
then $X$ will almost be contained in the range of $UQ_s U^*$ for
sufficiently large $s$ and hence a tiny (norm) perturbation of a
sufficiently large $UQ_s U^*$ will actually be the identity on $X$
(and still almost commute with $\mathfrak{F}$ and almost recover the
trace $\tau$).
\end{proof}

As mentioned above, the set $\UTAafd$ is the one most relevant to 
the classification program.  It also seems to be the most difficult 
one to understand in general.  However, for certain C$^*$-algebras it is 
easily seen to be all of $\TA$. 

\begin{defn}
We will call a C$^*$-algebra $A$ {\em tracially type I} if for each 
finite subset $\mathfrak{F} \subset A$ and $\varepsilon > 0$ there 
exists a type I subalgebra $B \subset A$ with unit $e$ such that 
$\| [x,e] \| < \varepsilon$ for all $x \in \mathfrak{F}$, 
$e\mathfrak{F} e \subset^{\varepsilon} B$ and $\tau(e) > 1 - 
\varepsilon$ for every $\tau \in \TA$.  
\end{defn}

One could also adopt Lin's notions used to define tracial topological
rank.  However, for our purposes the simpler definition above is
sufficient.  It is not hard to verify that if $A$ is a C$^*$-algebra
with finite tracial topological rank (either \cite[Definition
3.1]{lin:tracialtopologicalrank} or \cite[Definition
3.4]{lin:tracialtopologicalrank}) then $A$ is tracially type I in the
sense described above.  Moreover, any ASH algebra (with or without
slow dimension growth) is evidently tracially type I.

\begin{prop}
\label{thm:traciallytypeI}
If $A$ is tracially type I then $\TA = \UTAafd$.
\end{prop}

\begin{proof} Let $A$ be tracially type I and $\tau \in \TA$ be arbitrary. 
If $\mathfrak{F} \subset A$ is arbitrary and $\varepsilon > 0$ is
given then we must produce a u.c.p. map $\phi : A \to M_n ({\mathbb
C})$ such that $\phi$ is $\varepsilon$-multiplicative on
$\mathfrak{F}$ and $| \tau(x) - tr_n \circ \phi (x) | < \varepsilon$
for all $x$ in the unit ball of $A$.

Choose $\delta > 0$ very small and apply the definition of tracially
type I to $\mathfrak{F}$ and $\delta$ to get a type I subalgebra $B
\subset A$ with unit $e$.  Note that for every positive $x$ in the
unit ball of $A$ we have $| \tau(x) - \tau(exe) | < \delta$.  Also,
the map $x \mapsto exe$ is almost multiplicative on $\mathfrak{F}$.

If $\pi_{\tau} : A \to B(L^2(A,\tau))$ is the GNS representation, then
the weak closure of $\pi_{\tau} (B)$ is a type I von Neumann algebra
with faithful, normal trace given by $\tau(\cdot)/ \tau(e)$.

{\noindent Case 1:} Assume that the weak closure of $\pi_{\tau} (B)$
is isomorphic to $$\prod_{i = 1}^{k} M_{n(i)} ({\mathbb C}) \otimes
L^{\infty}(X_i, \mu_i)$$ for some probability spaces $(X_i, \mu_i)$.
Then the weak closure is a (non-separable) AF algebra.  In particular,
we can find a finite dimensional subalgebra, $C$, of the weak closure
of $\pi_{\tau} (eAe)$ which almost contains (up to $\delta$) the set
$\pi_{\tau}(e\mathfrak{F}e)$ {\em in norm}.  We then get the desired
map $\phi$ by composing $x \mapsto \pi_{\tau} (exe)$ with a
$\tau(\cdot)/ \tau(e)$ preserving conditional expectation from the
weak closure of $\pi_{\tau} (eAe)$ to $C$ (actually, one must also
apply Lemma \ref{thm:tracepreservingembedding} to get from $C$ to a
matrix algebra).

{\noindent Case 2:} The only other possibility is that the weak
closure of $\pi_{\tau} (B)$ is isomorphic to $$\prod_{i = 1}^{\infty}
M_{n(i)} ({\mathbb C}) \otimes L^{\infty}(X_i, \mu_i),$$ which is no
longer an AF algebra.  However, we are saved by the fact that
$\tau(\cdot)/ \tau(e)$ is a normal trace and hence we have $$1 =
\sum\limits_{i = 1}^{\infty} \frac{\tau(e_i)}{\tau(e)},$$ where $e_i$
denotes the unit of $M_{n(i)} ({\mathbb C}) \otimes L^{\infty}(X_i,
\mu_i)$.  Putting $E_n = e_1 + \ldots + e_n$ we then get that
$\tau(E_n)/ \tau(e) \to 1$.  Moreover, since $\pi_{\tau}
(e\mathfrak{F}e)$ is almost contained in $\pi_{\tau} (B)$ it follows
that the map $x \mapsto E_n \pi_{\tau} (exe) E_n$ is still almost
multiplicative (in norm) on $\mathfrak{F}$.  Taking $n$ large enough,
we now complete the proof as in case 1.
\end{proof}

As usual,  groups provide some interesting and instructive 
examples. As we are sticking to unital algebras, we will only consider 
the discrete case.  

\begin{example} (Compare with \cite[Corollary 2.11]{bedos:hypertraces}.)
If $\Gamma$ is a discrete group then there is a sort of `all or
nothing' principle for the weakly approximately finite dimensional
traces on the {\em reduced} group C$^*$-algebra $C^*_r (\Gamma)$. More
precisely, we have $T(C^*_r (\Gamma)) = UT(C^*_r (\Gamma))_{{\rm
w-AFD}}$ if and only if $\Gamma$ is amenable and $T(C^*_r
(\Gamma))_{{\rm w-AFD}} = \emptyset$ if $\Gamma$ is not amenable.  To
see that this is the case, we first recall that if $\Gamma$ is
amenable then $C^*_r (\Gamma)$ is nuclear and hence every trace is in
$UT(C^*_r (\Gamma))_{{\rm w-AFD}}$.  On the other hand, if there
exists a trace $\alpha \in T(C^*_r (\Gamma))_{{\rm w-AFD}}$ then (in
the left regular representation) $\alpha$ extends to a hypertrace
$\varphi$ on $B(l^2(\Gamma))$.  A simple calculation then shows that
$\varphi$ defines a left invariant mean on $l^{\infty} (\Gamma)
\subset B(l^2(\Gamma))$ and hence $\Gamma$ is amenable. (See 
\cite[Proposition 2.12]{bedos:hypertraces} 
for analogues  of this to certain (twisted) crossed product
C$^*$-algebras.)
\end{example}

The previous example clarifies an observation of J. Rosenberg: If 
C$^*_r (\Gamma)$ is QD then $\Gamma$ is amenable.  Since $\TAafd 
\neq \emptyset$ for every QD C$^*$-algebra the example above shows 
that it takes much less than quasidiagonality to imply amenability 
for C$^*_r (\Gamma)$.

There are two natural traces on the full group C$^*$-algebra,
$C^*(\Gamma)$. Namely the one coming from the trivial one dimensional
representation (which is clearly in $UT(C^*(\Gamma))_{{\rm AFD}}$) and
the one which gives the left regular representation after applying the
GNS construction.  This latter trace will be denoted by $\tau$, but in
the following examples it is important to remember that we are {\em no
longer} talking about the {\em reduced} group C$^*$-algebra.

\begin{example}
If $\Gamma$ is a residually finite discrete group then the canonical
trace $\tau$ is always in $T(C^*(\Gamma))_{{\rm AFD}}$.  To see this,
we let $\Gamma_1 \trianglerighteq \Gamma_2 \trianglerighteq \ldots$ be
a descending sequence of normal subgroups each of which has finite
index in $\Gamma$ and such that their intersection is the neutral
element.  Let $\phi_n : C^* (\Gamma) \to B(l^2(\Gamma/\Gamma_n))$ be
the unitary representation induced by the left regular representation
of $\Gamma/\Gamma_n$.  Since $l^2(\Gamma/\Gamma_n)$ is a finite
dimensional Hilbert space it follows that $\tau \in
T(C^*(\Gamma))_{{\rm AFD}}$.  Also note that $\tau$ is {\em uniformly}
weakly approximately finite dimensional if and only if $\Gamma$ is
amenable, by Theorem \ref{thm:mainthmUTAwafd}.  It follows that the
sets $\UTAwafd$ and $\UTAafd$ need not be closed in the weak-$*$
topology, in general.
\end{example}

\begin{example}
If $\Gamma$ is a discrete group with Kazhdan's Property T then
$T(C^*(\Gamma))_{{\rm AFD}} = T(C^*(\Gamma))_{{\rm w-AFD}}$ and,
moreover, $\tau \in T(C^*(\Gamma))_{{\rm w-AFD}}$ if and only if
$\gamma$ is residually finite.  Both of these claims follow from
\cite[Proposition 2.3]{kirchberg:propertyTgroups}.  The basic idea is that
Property T implies that any u.c.p.\ map which is weakly almost
multiplicative can be approximated (in a probabalistic sense) by an
honest homomorphism.
\end{example}

\begin{example} 
\label{thm:freegroupexample}
Let ${\mathbb F}_{\infty}$ be the free group on (countably) infinitely
many generators.  Then $T(C^*({\mathbb F}_{\infty}))_{\rm w-AFD} =
T(C^*({\mathbb F}_{\infty}))_{\rm AFD} = [T(C^*({\mathbb
F}_{\infty}))_{\rm FD}]^{-}$, where $[T(C^*({\mathbb
F}_{\infty}))_{\rm FD}]^{-}$ denotes the weak-$*$ closure of the set of
traces on $C^*({\mathbb F}_{\infty})$ whose GNS representations are
finite dimensional.  This observation appears in Kirchberg's work on
Connes' embedding problem for II$_1$ factors (see \cite[Lemma
4.5]{kirchberg:invent}).  It is also shown in \cite{kirchberg:invent}
that Connes' embedding problem is equivalent to verifying the equation
$T(C^*({\mathbb F}_{\infty}))_{\rm w-AFD} = T(C^*({\mathbb
F}_{\infty}))$. Indeed, what Kirchberg observes in
\cite{kirchberg:invent} is that a trace $\tau$ on $C^*({\mathbb
F}_{\infty})$ is weakly approximately finite dimensional if and only
if there exists a $\tau$-preserving embedding $\pi_{\tau}(C^*({\mathbb
F}_{\infty}))^{\prime\prime} \hookrightarrow R^{\omega}$. (See also
\cite{hadwin}.)  From this we see that, in general, $\TAwafd \neq
\UTAwafd$ and $\TAafd \neq \UTAafd$. (Since the latter sets give
hyperfinite GNS representations, by Theorem \ref{thm:mainthmUTAwafd},
while there are plenty of examples of non-hyperfinite,
$R^{\omega}$-embeddable von-Neumann algebras.)  Compare with the
locally reflexive case; Theorem \ref{thm:locallyreflexive}.
\end{example}

\section{The Nuclear Case}

In this section we discuss the connection between approximately finite
dimensional traces and Elliott's conjecture that simple, separable,
nuclear C$^*$-algebras are classified by their K-theoretic and tracial
data. Our first goal is to show that a necessary condition for
Elliott's conjecture to hold is that every unital, nuclear, QD
C$^*$-algebra $A$ satisfies the equation $\TA = \UTAafd$ (even the
non-simple ones!).  Then, at the end of the section, we will point to
a few references which lead us to believe that knowing an equation like
$\TA = \UTAafd$ may eventually be part of the sufficiency as well.

Our first lemma states that when dealing with many classes of 
C$^*$-algebras, it often suffices to restrict attention to 
 Popa algebras.

\begin{lem}
\label{thm:lemmasection4}

If $\mathfrak{C}$ is a collection of C$^*$-algebras which contains
$\mathbb C$ and is closed under $i)$ increasing unions (i.e.\
inductive limits with injective connecting maps), $ii)$ quasidiagonal,
semi-split extensions (i.e.\ if $0 \to I \to E \to B \to 0$ is a
semi-split (cf.\ \cite{blackadar:book}), short exact sequence, $I$
contains an approximate unit {\em of projections} which is
quasicentral in $E$ and both $I, B \in \mathfrak{C}$ then $E \in
\mathfrak{C}$) and $iii)$ tensoring with finite dimensional matrix
algebras then the following are equivalent:

\begin{enumerate}
\item $\TA = \TAafd$ (resp.\ $\TA = \UTAafd$) for every QD $A \in
\mathfrak{C}$.

\item $\TA = \TAafd$ (resp.\ $\TA = \UTAafd$) for every residually
finite dimensional $A \in \mathfrak{C}$.
 
\item $\TA = \TAafd$ (resp.\ $\TA = \UTAafd$) for every Popa algebra
$A \in \mathfrak{C}$.
\end{enumerate}

If, moreover, the class $\mathfrak{C}$ is closed under tensor products
with (non-unital) abelian algebras (it actually suffices to know $A
\in \mathfrak{C} \Longrightarrow A\otimes C_0 ((0,1]) \in
\mathfrak{C}$) then the following are equivalent:

\begin{enumerate}
\item[(4)] $\TA = \TAwafd$ (resp.\ $\TA = \UTAwafd$) for every $A \in
\mathfrak{C}$.

\item[(5)] $\TA = \TAwafd$ (resp.\ $\TA = \UTAwafd$) for every
QD $A \in \mathfrak{C}$.

\item[(6)] $\TA = \TAwafd$ (resp.\ $\TA = \UTAwafd$) for every residually
finite dimensional $A \in \mathfrak{C}$.
 
\item[(7)] $\TA = \TAwafd$ (resp.\ $\TA = \UTAwafd$) for every Popa algebra
$A \in \mathfrak{C}$.
\end{enumerate}
\end{lem}

\begin{proof} We first prove the equivalence of (1) - (3) 
 and then indicate the changes necessary to prove the second part.
 The proofs are the same whether dealing with $\TAafd$ or $\UTAafd$
 and hence we just treat the uniformly approximately finite
 dimensional case.

(1) $\Longrightarrow$ (3) is immediate.  (3) $\Longrightarrow$ (2)
    follows from Theorem \ref{thm:basicconstruction}.

(2) $\Longrightarrow$ (1).  Let $A \in \mathfrak{C}$, $\tau \in \TA$
    and ${\mathfrak{F}} \subset A$ be an arbitrary finite set.  Since
    $A$ is QD we can find a sequence of u.c.p. maps $\varphi_{n} : A
    \to M_{k(n)} ({\mathbb C})$ which are asymptotically
    multiplicative and asymptotically isometric (i.e.\ $\| a \| = \lim
    \| \varphi_{n} (a) \|$, for all $a \in A$).  Passing to a
    subsequence, if necessary, we may assume that $\varphi_{1}$ (and
    all the other $\varphi_n$'s) is as close to multiplicative on
    $\mathfrak{F}$ as we like.  Put $\Phi = \oplus_n \varphi_{n} : A
    \to \Pi_n M_{k(n)} ({\mathbb C})$, and let $E$ be the
    C$^*$-algebra generated by $\Phi(A)$.  Note that $\Phi : A \to E$
    is as close to multiplicative on $\mathfrak{F}$ as we like, by
    construction. 

Now observe that we have a semi-split, quasidiagonal, short exact
sequence: $$0 \to \oplus M_{k(n)} ({\mathbb C}) \to E + \oplus
M_{k(n)} ({\mathbb C}) \to A \to 0.$$ Since $\mathfrak{C}$ is closed
under all of the operations used, it follows that $E + \oplus M_{k(n)}
({\mathbb C}) \in {\mathfrak{C}}$ and it is clear that $ E + \oplus
M_{k(n)} ({\mathbb C})$ is residually finite dimensional.  Hence every
trace on $ E + \oplus M_{k(n)} ({\mathbb C})$ is uniformly
approximately finite dimensional by (2).  In particular, the trace
$\tau \in \TA$ defines a trace on $ E + \oplus M_{k(n)} ({\mathbb C})$
which is uniformly approximately finite dimensional.  Since the
splitting $\Phi : A \to E \subset E + \oplus M_{k(n)} ({\mathbb C})$
is almost multiplicative, it follows that we can construct a
u.c.p. map on $A$ (by composing maps on $E$ with $\Phi$) which is
almost multiplicative on $\mathfrak{F}$ and which approximately
recaptures the trace $\tau$.  This completes the proof of (2)
$\Longrightarrow$ (1).

For the equivalence of (4) - (7), we really only need to show the
implication (5) $\Longrightarrow$ (4) as the arguments above go
through without change for the other implications.  To prove this we
will need the following fact which follows from Theorems
\ref{thm:mainthmsection3} (part (8)) and \ref{thm:mainthmUTAwafd}
(part (5)): If $\pi:B \to A$ is a {\em surjective} $*$-homomorphism
and $\phi:A \to B$ is a c.p. splitting (i.e.\ $\pi \circ \phi = id_A$)
then for each $\tau \in \TA$ we have that $\tau \in \TAwafd$ (resp.\
$\tau \in \UTAwafd$) if and only if $\tau\circ\pi \in T(B)_{\rm
w-AFD}$ (resp.\ $\tau\circ\pi \in UT(B)_{\rm w-AFD}$).  (Thanks to
Lemma \ref{thm:nonunital} it does not matter whether or not the
splitting $A \to B$ is unital. Note also that the existence of a 
c.p.\ splitting is not necessary in the case of uniformly weakly 
approximately finite dimensional traces.)

With this observation in hand, (5) $\Longrightarrow$ (4) becomes very
simple.  Indeed, let $B$ be the unitization of the cone over $A$
(i.e.\ the unitization of $C_0 ((0,1]) \otimes A$).  Then $B$ is QD
(cf.\ \cite{dvv:QDhomotopy}) and belongs to $\mathfrak{C}$.  Moreover,
there is a natural surjective $*$-homomorphism $B \to A \oplus
{\mathbb C} \to A$.  A (non-unital) c.p. splitting for this quotient
map is given by $a \mapsto e \otimes a$, where $e \in C_0 ((0,1])$ is
any non-negative function such that $e(1) = 1$.  Hence if every trace
on $B$ is (uniformly) weakly approximately finite dimensional then
every trace on $A$ enjoys the same property. 
\end{proof}

The assumptions on the class $\mathfrak{C}$ may seem unusual, but note
that any one of the following classes of C$^*$-algebras is closed
under the operations needed in the lemma above: nuclear
C$^*$-algebras, exact C$^*$-algebras (cf.\ \cite[Section
7]{kirchberg:commutantsofunitaries}), real rank zero C$^*$-algebras
(cf.\ \cite[2.10, 3.1, 3.14]{brown-pedersen}, it is easy to prove that
if an extension is semi-split and quasidiagonal then every projection
in the quotient lifts to a projection in the middle algebra) - though
these algebras are not closed under tensoring with $C_0 ((0,1])$.

\begin{prop}
\label{thm:elliott}
If Elliott's conjecture holds for all nuclear, simple,
QD C$^*$-algebras with stable rank one and unperforated
K-theory then $\TA = \UTAafd$ for every nuclear,  QD
C$^*$-algebra $A$.
\end{prop}
 
\begin{proof} 
We will apply the previous lemma to the set $\mathfrak{C}$ of all
 nuclear C$^*$-algebras. We remark that extensions of nuclear
 C$^*$-algebras are again nuclear by \cite[Corollary
 3.3]{choi-effros:nuclearityandinjectivity}.

So assume that Elliott's Conjecture holds for all nuclear,
simple,  QD C$^*$-algebras with stable rank one and
unperforated K-theory and let $A$ be QD and
nuclear.  By the previous lemma, we may assume that $A$ is simple.
Let $B$ be a UHF algebra and note that every trace on $A$ extends (in
fact, uniquely) to a trace on $A \otimes B$. Hence it suffices to show
that $T(A\otimes B) = UT(A\otimes B)_{{\rm AFD}}$.  However, by 
\cite{rordam:tensorUHFI} and \cite{rordam:tensorUHFII},
$A\otimes B$ has stable rank one and unperforated K-theory.  Thus
$A\otimes B$ is classifiable.  But then, just as in the proof of 
Proposition \ref{thm:elliottpredictsASH}, it follows that $A\otimes B$
is an ASH algebra. Hence, by Proposition \ref{thm:traciallytypeI},
we see that $T(A\otimes B) = UT(A\otimes B)_{{\rm AFD}}$.
\end{proof}

The point of the above result is that verifying the equation $\TA =
\UTAafd$ for every nuclear, unital, QD C$^*$-algebra $A$ is a
necessary condition for Elliott's conjecture to hold.  However, we
also believe that knowing the equation $\TA = \UTAafd$ for a
particular $A$ may someday be part of the sufficient conditions for
classification. In \cite{lin:TAFclassification} Huaxin Lin proved that
certain tracially AF algebras are classifiable.  One of the key
technical tools he needed in the proof (see \cite[Lemmas 2.7 and
2.10]{lin:TAFclassification}) was a good understanding of the
approximation properties of traces on tracially AF algebras.  Hence we
would suggest that it is natural to try to classify algebras which
satisfy some form of tracial approximation property (something like
the equation $\TA = \UTAafd$) as a replacement for, say, assuming
tracially AF or ASH.  Indeed, the second part of \cite[Theorem
3.3]{popa:simpleQD} shows that a simple, QD C$^*$-algebra with real
rank zero and satisfying the equation $\TA = \UTAafd$ has a finite
dimensional approximation property which is vaguely similar to being
tracially AF (even without the assumption of nuclearity).  For
example, we have the following result. 

\begin{thm}
\label{thm:idontknow}
Let $A$ be a locally reflexive Popa algebra with real rank zero and
unique trace $\tau$. Then for every finite set ${\mathcal F} \subset A$
and $\varepsilon > 0$ there exists $n \in {\mathbb N}$, subalgebras
$Q_1, \ldots, Q_m \subset A$ each of which is isomorphic to $M_n
({\mathbb C})$ and with units $e_1, \ldots, e_m$ such that $\tau(e_1)
= \tau(e_2) = \cdots = \tau(e_m)$, $\| [e_i, x] \| < \varepsilon$ for
$1 \leq i \leq m$ and all $x \in {\mathcal F}$, $e_i {\mathcal F} e_i
\in^{\varepsilon} Q_i$ for $1 \leq i \leq m$ and, finally, $$ \|
\frac{1}{m\tau(e_1)} \sum_{k = 1}^{m} e_k - 1_A \|_{\tau,1} <
\varepsilon,$$ where $\| x \|_{\tau,1} = \tau(|x|)$.
\end{thm}

\begin{proof} 
By \cite[Theorem 3.3]{popa:simpleQD} it suffices to show that $\tau
\in \UTAafd$.  However, we will see later that the assumption of local
reflexivity implies that $\TAafd = \UTAafd$ (cf.\ Theorem
\ref{thm:locallyreflexive}) and hence we are done.
\end{proof}

Note that if one could somehow adapt Popa's techniques to further
arrange that $\tau(e_1) = 1/m$ then the approximation property above
would look very similar to a tracially AF algebra.

As mentioned in the proof above, $\TAafd = \UTAafd$ for every nuclear
C$^*$-algebra and hence deciding whether $\TA = \UTAafd$ is equivalent to
deciding if $\TA = \TAafd$ for every nuclear, QD C$^*$-algebra.  If
the hyperfinite II$_1$ factor is a QD C$^*$-algebra then this question
has an affirmative answer (see the appendix).

\section{The Exact Case}

The previous section attempts to argue that deciding whether or not
$\TA = \UTAafd$ for every {\em nuclear}, QD C$^*$-algebra is a natural
and important open question.  In this section we will show that this
equation fails in the category of {\em exact} C$^*$-algebras.  In
fact, it is possible to construct very nice exact C$^*$-algebras for
which $\TA \neq \TAwafd$.  Such algebras are not tracially type I (in
particular, do not have finite tracial topological rank in the sense
of \cite{lin:tracialtopologicalrank}) by Proposition
\ref{thm:traciallytypeI}.

We won't actually need condition (3) in the next proposition.
However, we feel that it may be of independent interest.  We wish to
thank G.\ Pisier for providing the first proof of (3)
$\Longrightarrow$ (4), though Ozawa later observed that exactness was
an unnecessary assumption (see Theorem \ref{thm:mainthmUTAwafd}).

\begin{prop} 
\label{prop:mainpropsection4}
Let $A$ be a unital, exact C$^*$-algebra and $\tau \in \TA$. Then the
following are equivalent:

\begin{enumerate}
\item $\tau \in \TAwafd$.

\item $\tau \in \UTAwafd$.

\item The GNS representation is a nuclear map into $\pi_{\tau}
(A)^{\prime\prime}$.  That is, there exist u.c.p. maps $\phi_n : A \to
M_{k(n)} ({\mathbb C})$ and $\psi_n : M_{k(n)} ({\mathbb C}) \to
\pi_{\tau} (A)^{\prime\prime}$ such that $\| \pi_{\tau} (a) - \psi_n
\circ \phi_n (a) \| \to 0$ for all $a \in A$.

\item $\pi_{\tau} (A)^{\prime\prime}$ is hyperfinite. 
\end{enumerate}
\end{prop}

\begin{proof} We will see in the next section that (1) and (2) are 
equivalent since Kirchberg has shown that exactness implies local
reflexivity.  Hence by Theorem \ref{thm:mainthmUTAwafd} we see that
(1), (2) and (4) are equivalent.  Since nuclearity obviously implies
weak nuclearity we also get the implication (3) $\Longrightarrow$ (4)
from Theorem \ref{thm:mainthmUTAwafd}.  Hence we are only left to
prove (1) $\Longrightarrow$ (3).

We first recall one of Kirchberg's characterizations of exactness: $A$
is exact if and only if every u.c.p. map $A \to D$, where $D$ is an
injective C$^*$-algebra, is nuclear.  From this fact, we see that it
suffices to show that the GNS representation $\pi_{\tau} : A \to
\pi_{\tau} (A)^{\prime\prime}$ factors through an injective
C$^*$-algebra.  More precisely, it suffices to show that there exist
u.c.p. maps $\Phi : A \to l^{\infty} (R)$ and $\Psi : l^{\infty} (R)
\to \pi_{\tau} (A)^{\prime\prime}$ such that $\pi_{\tau} = \Psi \circ
\Phi$, where $R$ denotes the hyperfinite II$_1$ factor.  Indeed, if we
can do this, then the map $\Phi$ is nuclear, by Kirchberg's
characterization of exactness, and hence $\pi_{\tau} = \Psi \circ
\Phi$ is also nuclear.

So assume that $\tau \in \TAwafd$.  Then there exists a u.c.p. map
$\Phi : A \to l^{\infty} (R)$ such that $\sigma \circ \Phi$ is a
$*$-homomorphism with $\tau = \tau_{\omega} \circ \sigma \circ \Phi$,
where $\sigma : l^{\infty} (R) \to R^{\omega}$ is the quotient mapping
and $\tau_{\omega}$ is the unique trace on the ultrapower of the
hyperfinite II$_1$ factor.  Since the weak closure of $\sigma \circ
\Phi (A)$ is canonically isomorphic to $\pi_{\tau}
(A)^{\prime\prime}$, and there is a conditional expectation $E:
R^{\omega} \to \pi_{\tau} (A)^{\prime\prime}$, we get the desired map
by defining $\Psi = E \circ \sigma$.
\end{proof}

We now come to the main result of this section.  

\begin{thm}
\label{thm:exoticpopaalgebra}
There exists an exact, Popa algebra, $A$, with real rank zero, stable
rank one, UCT, Blackadar's fundamental comparison property (i.e.\ if
$p,q \in A$ are projections such that $\tau(q) < \tau(p)$ for all
$\tau \in T(A)$ then $q$ is equivalent to a subprojection of $p$),
unperforated K-theory, Riesz interpolation property and which is
approximately divisible and an increasing union of residually finite
dimensional subalgebras such that $\TA \neq \TAwafd$.
\end{thm}

\begin{proof} We first claim that it suffices to construct a 
C$^*$-algebra $C$ which is residually finite dimensional, exact, real
rank zero, satisfies the UCT and such that $T(C) \neq T(C)_{\rm
w-AFD}$.  Indeed, if we can find such a $C$ then by applying Theorem
\ref{thm:basicconstruction} to the class of all exact, real rank zero
C$^*$-algebras which satisfy the UCT we can find an exact Popa algebra
with real rank zero, UCT and such that $\TA \neq \TAwafd$.  (See
Remark \ref{thm:UCTremark} for the UCT assertion and
\cite{brown-pedersen} for a proof that matrices over a real rank zero
algebra also have real rank zero.)  Then replacing $A$ with $A\otimes
\mathcal{U}$, where $\mathcal{U}$ is some UHF algebra, will be the
desired example since this operation preserves Popa's property,
exactness, real rank zero, UCT and picks up stable rank one, Riesz
interpolation (cf.\ \cite[Corollary 3.15]{blackadar-kumjian-rordam}),
Blackadar's fundamental comparison property and hence unperforated
K-theory (cf.\ \cite{rordam:tensorUHFI}, \cite{rordam:tensorUHFII}).
Moreover, it is clear that this example will be an inductive limit of
residually finite dimensional subalgebras and be approximately
divisible in the sense of \cite{blackadar-kumjian-rordam}.

The construction of the desired residually finite dimensional
C$^*$-algebra is a consequence of another one of Kirchberg's
characterizations of exactness \cite[Theorem
1.3]{kirchberg:commutantsofunitaries}: A separable C$^*$-algebra $A$
is exact if and only if there exists a subalgebra $B$ of the CAR
algebra, $M_{2^{\infty}}$, and an AF ideal $J \subset B$ such that $A
\cong B/J$.  We remark if $A$ is exact and $0 \to J \to B \to A \to 0$
is the short exact sequence given by Kirchberg's theorem, then this
sequence is automatically semi-split (i.e.\ there exists a
c.p. splitting $A \to B$) since $B$ is exact and $J$ is nuclear (cf.\
the bottom of page 41 in \cite{kirchberg:commutantsofunitaries}).

Since $C^*_r ({\mathbb F}_2)$ has a unique trace, it follows from
\cite[Theorem 7.2]{rordam:tensorUHFII} that $C^*_r ({\mathbb
F}_2)\otimes M_{2^{\infty}}$ is exact and has real rank zero.  By
Kirchberg's characterization, we can find an exact, QD C$^*$-algebra
$B$ with an AF ideal $J \subset B$ such that $B/J \cong C^*_r
({\mathbb F}_2)\otimes M_{2^{\infty}}$.  Moreover, the short exact
sequence $0 \to J \to B \to C^*_r ({\mathbb F}_2)\otimes
M_{2^{\infty}} \to 0$ is semi-split.

Since $J$ is AF, it follows from \cite[Theorem 3.14 and Corollary
3.16]{brown-pedersen} that $B$ also has real rank zero.  Note that
from part (4) of Proposition \ref{prop:mainpropsection4} it follows
that $T(B) \neq T(B)_{{\rm w-AFD}}$. Finally, since $C^*_r ({\mathbb
F}_2)$ satisfies the UCT, it follows from the `two out of three
principle' (cf.\ Proposition \ref{thm:twooutofthreeprinciple}) that
$B$ also satisfies the UCT.  $B$ is almost the desired algebra; we
only have to replace $B$ with something residually finite dimensional
($B$ is only QD).

From the proof of (2) $\Longrightarrow$ (1) in Lemma
\ref{thm:lemmasection4} we can use $B$ to construct a residually
finite dimensional C$^*$-algebra $C$ such that $C$ is exact, real rank
zero, satisfies the UCT and such that $T(C) \neq T(C)_{{\rm w-AFD}}$. 
\end{proof}

\begin{cor}
\label{thm:notTAF}
There exists an exact, Popa algebra, $A$, with real rank zero, stable
rank one, UCT, unperforated K-theory, Riesz interpolation and
Blackadar's fundamental comparison property which is approximately
divisible and an increasing union of residually finite dimensional
subalgebras but such that $A$ is not tracially AF.
\end{cor}

\begin{proof}
Since $\TA = \UTAafd$ for every tracially type I algebra (see 
Proposition \ref{thm:traciallytypeI}), the example in 
the previous theorem can't be tracially AF.
\end{proof}

It may seem unusual to mention in the previous two results that $A$ is
an inductive limit of residually finite dimensional subalgebras.  Our
main reason for pointing out this fact is that some of Lin's recent
structural work on the class of tracially AF C$^*$-algebras relies
heavily on a theorem of Blackadar and Kirchberg stating that every
simple, nuclear, QD C$^*$-algebra is an inductive limit of such
subalgebras.  In fact, for some of Lin's structural work, this is the
only place that nuclearity is used (i.e.\ his results hold more
generally if one replaces the assumption of nuclearity by the
assumption of an inductive limit decomposition by residually finite
dimensional subalgebras).

In \cite{bedos:hypertraces} B\'{e}dos asked whether or not every
separable, unital hypertracial C$^*$-algebra is nuclear (see
\cite[Section 3]{bedos:hypertraces} - in the language of the present
paper, a C$^*$-algebra is hypertracial if every quotient has at least
one weakly approximately finite dimensional trace).  It is easy to see
that every simple, unital, QD C$^*$-algebra is hypertracial and hence
Dadarlat's examples of non-nuclear Popa algebras provide
counterexamples to this question
\cite{dadarlat:nonnuclearsubalgebras}.  Theorem
\ref{thm:exoticpopaalgebra} above provides further examples.  Indeed,
for every non-hyperfinite II$_1$ factor $M$ which contains a weakly
dense exact C$^*$-subalgebra the proof of Theorem
\ref{thm:exoticpopaalgebra} shows that we can construct an exact Popa
algebra with stable rank one, Blackadar's comparison property (hence
unperforated K-theory), Riesz property, approximate divisibility and
which is an increasing union of residually finite dimensional
subalgebras but which is not nuclear since it will have
$M\bar{\otimes} R$ as the weak closure of some GNS representation.

\section{The Locally Reflexive Case}

The main result of this section is: 

\begin{thm} 
\label{thm:locallyreflexive}
If $A$ is locally reflexive then $\TAwafd = \UTAwafd$ and $\TAafd =
\UTAafd$.
\end{thm}

\begin{proof}
We only give the proof of $\TAafd = \UTAafd$ as it will be clear that
essentially the same proof gives the other equality. 

So let $\tau \in \TAafd$ be arbitrary.  Evidently it suffices to prove
that if $\mathfrak{F} \subset A$ is an arbitrary finite set and
$\varepsilon > 0$ then there exists a u.c.p. map $\varphi : A \to B$,
where $B$ is a finite dimensional C$^*$-algebra, such that $\|
\varphi(xy) - \varphi(x)\varphi(y) \| < \varepsilon$, for all $x,y \in
\mathfrak{F}$ and such that there exists a trace $\gamma$ on $B$ such
that $\| \tau - \gamma\circ\varphi \|_{A^*} < \varepsilon$.  

In order to do this, we will show that for each finite dimensional
operator system $X \subset A^{**}$ containing both the set $\mathfrak{F}$ 
and $\{ ab : a,b \in \mathfrak{F} \}$,
there exists a sequence of normal, u.c.p. maps $\psi_n : A^{**} \to
M_{s(n)} ({\mathbb C})$ such that $tr_{s(n)} \circ \psi_n (x) \to
\tau^{**} (x)$, for all $x \in X$ and {\em each} $\psi_n$ is
$\varepsilon$-multiplicative on $\mathfrak{F}$.  If we are able to do
this then one can construct a net of normal, u.c.p. maps
$\varphi_{\lambda} : A^{**} \to M_{k(\lambda)} ({\mathbb C})$ with the
property that $tr_{k(\lambda)} \circ \varphi_{\lambda} \in A^*$ for
all $\lambda$, each $\varphi_{\lambda}$ is
$\varepsilon$-multiplicative on $\mathfrak{F}$ and (here is the key)
$tr_{k(\lambda)} \circ \varphi_{\lambda} \to \tau^{**}$ in the weak
topology coming from $A^{**}$.  Hence, by the Hahn-Banach theorem,
$\tau^{**}$ belongs to the {\em norm} closure of the convex hull of
$\{ tr_{k(\lambda)} \circ \varphi_{\lambda} \} \subset A^*$.  Then one
would be able to choose a finite set $\lambda_1, \ldots, \lambda_p$
and positive real numbers $\theta_1, \ldots, \theta_p$ such that $\sum
\theta_i = 1$ and $\| \tau^{**} - \sum \theta_i tr_{k(\lambda_i)}
\circ \varphi_{\lambda_i} \|_{A^*} < \varepsilon$.  Finally one
defines $B = M_{k(\lambda_1)}({\mathbb C}) \oplus \cdots \oplus
M_{k(\lambda_p)} ({\mathbb C})$, $\varphi = \varphi_{\lambda_1} \oplus
\cdots \oplus \varphi_{\lambda_p}$ and $\gamma = \sum \theta_i
tr_{k(\lambda_i)}$.  As we have arranged that $\varphi_{\lambda}$ is
$\varepsilon$-multiplicative on $\mathfrak{F}$ {\em for every}
$\lambda$, it is clear that $\varphi$ will also be close to
multiplicative on $\mathfrak{F}$.  

So let $X \subset A^{**}$ be any finite dimensional operator system
containing the sets $\mathfrak{F}$ and $\{ ab : a,b \in \mathfrak{F}
\}$.  Since $\tau \in \TAafd$ we can find a sequence of u.c.p. maps
$\varphi_m : A \to M_{k(m)} ({\mathbb C})$ which are asymptotically
multiplicative (in norm) and which recapture $\tau$ (as a weak-$*$
limit) after composing with the traces on $M_{k(m)} ({\mathbb C})$.
Note that by passing to a subsequence, if necessary, we may further
assume that each $\varphi_m$ is as close to multiplicative as one
likes on the set $\mathfrak{F}$. Since $A$ is locally reflexive, we
can find a net of u.c.p. maps $\alpha_i : X \to A$ such that $\alpha_i
(x) \to x$ in the weak-$*$ topology (coming from $A^*$) for all $x \in
X$.  Another standard Hahn-Banach argument allows us to extract a
subnet $\beta_t : X \to A$ from the convex hull of the $\alpha_i$'s
such that $\beta_t (x) \to x$ in the weak-$*$ topology, for all $x \in
X$, and $\| a - \beta_t (a) \| \to 0$ for all $a \in \mathfrak{F} \cup
\{ ab : a,b \in \mathfrak{F} \}$.  Again, passing to a subnet, we may
assume that $\| \beta_t (a) - a \| < \varepsilon$ for all $a \in
\mathfrak{F} \cup \{ ab : a,b \in \mathfrak{F} \}$ and for all $t$.
In particular, this implies that $\varphi_m \circ \beta_t$ is nearly
multiplicative on $\mathfrak{F}$ for all $t,m$.

We almost have the desired maps $\psi_n$. Since $X$ is finite
dimensional we can choose a linear basis $\{ x_1, \ldots, x_q \}$.
For each $n \in {\mathbb N}$ first choose $t(n)$ such that $| \tau^{**}
(x_i) - \tau(\beta_{t(n)}(x_i)) | < 1/n$ for $1 \leq i \leq q$.  Then
choose $m_n$ such that $| \tau(\beta_{t(n)}(x_i)) - tr_{k(m_n)} \circ
\varphi_{m_n} (\beta_{t(n)}(x_i)) | < 1/n$ for $1 \leq i \leq q$.  Then
defining $\tilde{\psi}_n = \varphi_{m_n} \circ \beta_{t(n)} : X \to
M_{k(m_n)} ({\mathbb C})$ we have that $tr_{k(m_n)} \circ
\tilde{\psi}_n (x) \to \tau^{**}(x)$ for all $x \in X$.

By Arveson's Extension Theorem, we may assume that each
$\tilde{\psi}_n$ is actually defined on all of $A^{**}$.  The only
problem is that we can't be sure that the Arveson extensions are
normal on $A^{**}$.  However a tiny perturbation of the
$\tilde{\psi}_n$ will yield normal maps $\psi_n$ with all the right
properties. (Use the fact that c.p. maps to matrix algebras are
nothing but positive linear functionals on matrices over the given
algebra and that the set of normal linear functionals on a von Neumann
algebra is dense in the dual space.)
\end{proof}

The following corollary generalizes \cite[Theorem
7.5]{kirchberg:invent} and the discrete case of \cite[Proposition
7.1]{kirchberg:invent}. 

\begin{cor} 
Let $\Gamma$ be a discrete group such that every finitely generated
subgroup has Kirchberg's factorization property (cf.\
\cite{kirchberg:propertyTgroups}: for example, if $\Gamma$ is an
increasing union of residually finite subgroups).  Then $\Gamma$ is
amenable if and only if $C^*(\Gamma)$ is locally reflexive.
\end{cor}

\begin{proof} 
According to \cite[Lemma 4.1]{kirchberg:propertyTgroups} the
factorization property for a discrete group $\Gamma$ is equivalent to
knowing that the canonical trace on $C^* (\Gamma)$ (which gives the
left regular representation after applying GNS) is weakly
approximately finite dimensional. (This is essentially a special case
of the equivalence of (1) and (5) in Theorem
\ref{thm:mainthmsection3}.) Thus our assumptions imply that the
canonical trace on $C^* (\Gamma)$ is weakly approximately finite
dimensional (cf.\ \cite[Lemma 7.3, part ($v$)]{kirchberg:invent}).

Now the proof is simple.  The `only if' part is immediate since $C^*
(\Gamma) = C^*_r (\Gamma)$ is nuclear whenever $\Gamma$ is amenable.
The opposite direction follows from Theorem \ref{thm:locallyreflexive}
since injectivity of the von Neumann algebra generated by a discrete
group in it's left regular representation implies amenability of the
group.
\end{proof}

\section{The WEP  Case}

Here we observe that if $A$ has the WEP of E.C. Lance then $\TA =
\TAwafd$.  As an application it follows that {\em none} of the known
examples of non-quasidiagonal C$^*$-algebras can be embed into a
nuclear, {\em stably finite} C$^*$-algebra (even though many are exact
and stably finite). This should be contrasted with Kirchberg's
embedding theorem which states that every exact C$^*$-algebra embeds
into the Cuntz algebra $\mathcal{O}_2$.

Recall that a C$^*$-algebra $A$ is said to have the {\em weak
expectation property} (WEP) if for every {\em faithful}, nondegenerate
representation $\pi : A \to B(H)$ there exists a u.c.p. map $\Phi :
B(H) \to \pi(A)^{\prime\prime}$ such that $\Phi(\pi(a)) = \pi(a)$, for
all $a \in A$.  This is a large class of C$^*$-algebras (including
every injective C$^*$-algebra and every nuclear C$^*$-algebra).  In
fact, \cite[Corollary 3.5]{kirchberg:invent} states that every
(separable) C$^*$-algebra is contained in a (separable) simple
C$^*$-algebra with the WEP.

Kirchberg studied this class of algebras in \cite{kirchberg:invent}
and proved the following result: $A$ has the WEP if and only if there
is a unique C$^*$-norm on $A \odot C^* ({\mathbb F}_{\infty})$, where
$C^* ({\mathbb F}_{\infty})$ denotes the full group C$^*$-algebra of
the free group on infinitely many generators.

\begin{prop}
\label{thm:WEPcase}
If $A$  has the WEP then $\TA = \TAwafd$.
\end{prop}

\begin{proof} 
 Let $\tau \in \TA$ and $\pi_{\tau}$ be the associated GNS
 representation.  Let $\rho : C^* ({\mathbb F}_{\infty}) \to
 \pi_{\tau}(A)^{\prime}$ be any *-homomorphism with weakly dense
 range.  By Kirchberg's tensorial characterization of the WEP we get a
 *-homomorphism $\pi_{\tau} \otimes \rho : A
 \otimes C^* ({\mathbb F}_{\infty}) \to
 C^*(\pi_{\tau}(A),\pi_{\tau}(A)^{\prime})$.  As in the proof of (5)
 $\Longrightarrow$ (6) from Theorem \ref{thm:mainthmsection3} it
 follows that $\tau$ is weakly approximately finite dimensional.
\end{proof}

Our next three corollaries are  related to some observations 
of Bed{\o}s (see, for example, \cite[Corollary 2.11]{bedos:hypertraces}). 

\begin{cor}
\label{thm:nohomomorphisms}
Let $A$ be such that $\TAwafd = \emptyset$ and $B$ be a C$^*$-algebra 
with the WEP and at least one tracial state.  Then there is no 
unital $*$-homomorphism $A \to B$.
\end{cor}

We already observed that the set of weakly approximately finite
dimensional traces is empty when $A = C^*_r (\Gamma)$ for a
non-amenable, discrete group $\Gamma$ and hence the corollary above
covers these examples.  These are the standard examples of stably
finite, non-QD C$^*$-algebras.  However, S. Wassermann has produced
other examples of stably finite, non-QD C$^*$-algebras in
\cite{wassermann:annals}, \cite{wassermann:nonQDstablyfinite} (though,
to the best of our knowledge, it is not known if any of these are
exact).  A careful inspection of Wassermann's proofs shows that he
proves much more than non-quasidiagonality; he actually shows that his
examples admit no weakly approximately finite dimensional traces and
hence are also covered by the corollary above.  At the moment the only
other known examples of stably finite, non-QD C$^*$-algebras are
things which contain one of the examples described above -- e.g.\ many
reduced free products with respect to tracial states contain reduced
free group C$^*$-algebras.  (Caution: There are plenty of stably
finite C$^*$-algebras for which it is not known if they are QD, but
the examples above are the only ones which are known to not be QD.)

\begin{cor}
None of the known examples (i.e.\ the examples described above) of
stably finite, non-QD C$^*$-algebras can be embed into a finite,
hyperfinite von Neumann algebra.
\end{cor}

Since (unital) stably finite, exact C$^*$-algebras always have at 
least one tracial state we have: 

\begin{cor}
None of the known examples of stably finite, non-QD C$^*$-algebras 
can be embed into a stably finite, nuclear C$^*$-algebra.
\end{cor}

The corollary above is well known to the free probability community in
the case of $C^*_r ({\mathbb F}_n)$ (or many other free product
groups).  Indeed, if $C^*_r ({\mathbb F}_n) \subset B$ and $B$ has a
trace $\tau$ then the restriction of $\tau$ to $C^*_r ({\mathbb F}_n)$
is the canonical trace (by uniqueness) and hence the weak closure of
$C^*_r ({\mathbb F}_n)$ in the GNS representation of $B$ with respect
to $\tau$ will not be hyperfinite and hence the weak closure of $B$
can't be hyperfinite.

Note, of course, that not only are embeddings impossible, but there
are no non-zero homomorphisms from any of the standard examples of
stably finite, non-QD C$^*$-algebras to a finite, hyperfinite
von-Neumann algebra (or stably finite, nuclear C$^*$-algebra).

The author's interest in embedding questions of the above type stems
from a question of Blackadar and Kirchberg which asks whether every
nuclear, stably finite C$^*$-algebra is QD (cf.\
\cite{blackadar-kirchberg}).  One consequence of the results above is
that it is impossible to give counterexamples to this question by
embedding one of the known non-QD C$^*$-algebras into a nuclear,
stably finite C$^*$-algebra.  Hence constructing a counterexample (if
one exists) will require completely new ideas.

Our final result of this section generalizes \cite[Theorem
6.8]{bekka}.  If $\pi : G \to B(H)$ is a unitary representation of a
locally compact group $G$ then Bekka defines the representation to be
amenable if there exists a state $\varphi$ on $B(H)$ such that
$\varphi(\pi(g) T \pi(g^{-1})) = \varphi(T)$ for all $T \in B(H)$ and
for all $g \in G$.  If $A = C^*( \{ \pi(g) : g \in G \} ) = span \{
\pi(g) : g \in G \}^{-}$ (norm closure, since $G$ is a group) and $A$
falls in the centralizer of some state on $B(H)$ then clearly the
representation is amenable.  However, the converse is also true.
Namely, if $\pi$ is amenable and $\varphi$ is a state such that
$\varphi(\pi(g) T \pi(g^{-1})) = \varphi(T)$ then we get
$$\varphi(\pi(g) T) = \varphi(\pi(g^{-1})[\pi(g) T] \pi(g)) =
\varphi(T\pi(g))$$ for all $T \in B(H)$.  Since $A$ is the norm
closure of the span of $\pi(G)$ it follows that $A$ also falls in the
centralizer of $\varphi$.  In summary: A representation $\pi : G \to
B(H)$ is amenable in the sense of Bekka if and only if $\TAwafd \neq
\emptyset$, where $A = C^*( \{ \pi(g) : g \in G \} )$.  In light of
Proposition \ref{thm:WEPcase}, the following corollary in now
immediate.

\begin{cor}
Let $\pi : G \to B(H)$ be a strongly continuous, unitary representation 
of a locally 
compact group $G$ and $A = C^*( \{ \pi(g) : g \in G \} )$.  If 
$A$ has the WEP then $\pi$ is an amenable representation if and 
only if $A$ has a tracial state.
\end{cor}

\section{II$_1$ Factor Representations of Popa Algebras}

The main motivation for Popa's work in \cite{popa:simpleQD} was to try
to understand the relationship between quasidiagonality and
nuclearity.  Indeed, in \cite{popa:simpleQD} Popa asked whether every
Popa algebra with unique trace was nuclear.  More generally, he also
asked whether $R$ was the only II$_1$ factor that could be realized as
the weak closure of the GNS representation of a Popa algebra.  There
was serious speculation that both questions should be true.  However
counterexamples to Popa's first question were constructed by Dadarlat
in \cite{dadarlat:nonnuclearsubalgebras}. Interestingly enough,
though, in \cite{dadarlat:nonnuclearsubalgebras} Dadarlat constructs
nonnuclear tracially AF algebras and hence all of their II$_1$ factor
representations are hyperfinite (even though the algebras are not
nuclear).  Further support for Popa's second question was also
provided in \cite[Remark 3.4.2]{popa:simpleQD} where Popa proved that
if a factorial trace is in $\UTAafd$ then it gives the hyperfinite
II$_1$ factor (compare with Theorem \ref{thm:mainthmUTAwafd}).
However, the results of the previous sections also answer Popa's
second question.

\begin{prop}
There exists an exact, Popa algebra, $A$, with real rank zero, stable
rank one, UCT, unperforated K-theory, Riesz interpolation and
Blackadar's fundamental comparison property which is approximately
divisible and an increasing union of residually finite dimensional
subalgebras but such that $A$ has non-hyperfinite II$_1$ factor
representations.
\end{prop}

\begin{proof}
Let $A$ be the example constructed in Theorem
\ref{thm:exoticpopaalgebra}.  Then any extreme point of $\TA$ which
does not belong to $\TAwafd$ will give a non-hyperfinite II$_1$ factor
representation.
\end{proof}

Note that if $A$ is any QD, locally reflexive C$^*$-algebra then $A$
has at least one hyperfinite II$_1$ (or finite dimensional) factor
representation. (Use the fact that $\TAwafd$ is a face in $\TA$, hence
contains an extreme point of $\TA$  and, finally, Theorem 
\ref{thm:locallyreflexive}.) In particular, the following result
shows that Popa's II$_1$ factor representation question has an
affirmative answer in the locally reflexive, unique trace case.

\begin{cor}
\label{thm:locallyreflexiveuniquetrace}
If $A$ is a locally reflexive, QD C$^*$-algebra with unique trace $\tau$
then either $\pi_{\tau} (A)^{\prime\prime} \cong R$ or 
$\pi_{\tau} (A)^{\prime\prime}$ is a finite dimensional matrix algebra. 
\end{cor}

\begin{rem}
One of the deficiencies of Theorem \ref{thm:basicconstruction} is that
it will usually produce a tracial space with an infinite number of
extreme points.  Since we know that the unique trace case always
produces $R$ (for infinite dimensional, locally reflexive Popa
algebras), a natural question then becomes: Is there an exact Popa
algebra such that $\TA$ has a finite number of extreme points, some of
which do not give hyperfinite GNS representations?

It is also natural to wonder whether local reflexivity is really needed 
in the corollary above.  That is, if a Popa algebra has unique trace 
then must one get $R$ in the GNS representation?
\end{rem}

Dropping the assumption of exactness, it is now easy to construct a 
sort of universal Popa algebra which realizes every McDuff factor as 
a GNS representation.  (Recall that a II$_1$ factor is called McDuff 
if it is isomorphic to something of the form $M\bar{\otimes}R$.)

\begin{thm} 
\label{thm:arbitraryMcDuff}
There exists a Popa algebra $A$ with the property that for each McDuff
factor, $M$, there exists a trace $\tau_M \in \TA$ such that
$\pi_{\tau} (A)^{\prime\prime} \cong M$.
\end{thm}

\begin{proof}
Since every II$_1$ factor arises as the weak closure of a GNS
representation of $C^* ({\mathbb F}_{\infty})$, and $C^* ({\mathbb
F}_{\infty})$ is residually finite dimensional, this theorem follows
from Theorem \ref{thm:basicconstruction}.
\end{proof}

\begin{cor}
If $M \subset B(L^2(M))$ is a McDuff factor then there exists a weakly
dense C$^*$-subalgebra $A \subset M$ such that:

\begin{enumerate}
\item For each finite subset $\mathfrak{F} \subset A$ and $\varepsilon
> 0$ there exists a finite dimensional subalgebra $B \subset A$ with
unit $e$ such that $\| E_B(exe) - (x - e^{\perp} x e^{\perp}) \|_2 <
\varepsilon \| e \|_2$, where $E_B : eMe \to B$ is a trace
preserving conditional expectation.

\item There exist finite rank projections, $P_1 \leq P_2 \leq \ldots$,
such that $\| [P_n, a] \| \to 0$ for all $a \in A$. 

\item There exists a state $\phi$ on $B(L^2(M))$ such that $A \subset
B(L^2(M))_{\phi} = \{ T \in B(L^2(M)): \phi(TS) = \phi(ST), S \in
B(L^2(M))\}.$

\item There exists a sequence of II$_1$ factors, $R_n \subset
B(L^2(M))$, such that $R_n \cong R$ for all $n$ and for each $a \in A$
we can find $x_n \in R_n$ such that $\| a - x_n \| \to 0$.

\item There exists a sequence of normal, u.c.p. maps $\varphi_n : M
\to M_{k(n)} ({\mathbb C})$ such that $\| \varphi_n (ab) - \varphi_n
(a) \varphi_n (b) \|_2 \to 0$ for all $a,b \in A$.

\item There exists a completely positively liftable u.c.p. map $\Phi :
M \to R^{\omega}$ such that $\Phi|_A$ is a $*$-monomorphism.

\end{enumerate}
\end{cor}

\begin{proof}
Let $A$ be the universal Popa algebra from the previous theorem.
Evidently $A$ satisfies (1) above since it satisfies the stronger norm
approximation property. This also immediately gives (5).  Since $A$ is
QD we get (2) and the fact that $\TAafd$ is not empty.  Hence, from
Theorem \ref{thm:mainthmsection3}, we get (3).  Since $A$ also embeds
into $R$, we get (4) from Voiculescu's Theorem.  Finally, note that
(6) follows from (5).
\end{proof}

When we first discovered part (6) in the corollary above, we thought
that it would be useful in showing that every II$_1$ factor embeds
into $R^{\omega}$.  However, if the map $\Phi$ is normal then $M \cong
R$.

Another curious consequence of this work is that McDuff factors 
which are generated by
exact C$^*$-algebras always have `norm microstates'  on a
dense subalgebra. 

\begin{thm}
If $M \subset B(L^2(M))$ is McDuff and contains a weakly dense, exact
C$^*$-subalgebra then there exists a weakly dense C$^*$-subalgebra $A
\subset M$ and finite dimensional matrix subalgebras $M_n \subset
B(L^2(M))$ such that for each $a \in A$ there exists a sequence $a_n
\in M_n$ such that $\| a - a_n \| \to 0$.  (Hence, for every
noncommutative polynomial $P$ in $k$ variables and finite set $\{
a^{(i)} \}_{i = 1}^{k} \subset A$ we have $\| P(a^{(1)}, \ldots,
a^{(k)}) - P(a^{(1)}_n, \ldots, a^{(k)})_n) \| \to 0$ as $n \to
\infty$.)
\end{thm}

\begin{proof}
Since every exact C$^*$-algebra is the quotient of an exact,
residually finite dimensional C$^*$-algebra (cf.\
\cite[Corollary 5.3]{brown:QDsurvey}), it follows from Theorem
\ref{thm:basicconstruction} that $M$ contains a weakly dense, exact
Popa algebra.  In particular, $M$ contains a weakly dense, exact, QD
C$^*$-algebra.  Since $M$ is a factor, it can't contain any nonzero
compact operators and hence the result now follows from
\cite{dadarlat:QDapproximation} (see also \cite{brown:herrero} for the
general case).
\end{proof}

Note that the preceding theorem covers many group von Neumann algebras
(e.g.\ $\Gamma = G_1 \times G_2$ where $G_1$ is discrete, amenable and
i.c.c. while $G_2$ is discrete, exact and i.c.c.). We are not,
however, claiming that this result implies $R^{\omega}$ embeddability
for such group von Neumann algebras.  Indeed, it is not at all clear
that the existence of norm microstates implies the existence of `weak'
microstates (in the sense of Voiculescu) since there does not appear
to be any way of understanding how the traces behave on the norm
approximations.

\section{Connes' Embedding Problem}
In this section we show that the techniques of this paper easily
yield a new characterization of those McDuff factors which are
embeddable into $R^{\omega}$.  An open problem of Connes asks whether
or not every (separable) II$_1$ factor embeds into $R^{\omega}$ and it
is easy to see that this is the case if and only if every McDuff
factor embeds into $R^{\omega}$.  Embeddable II$_1$ factors already
admit a number of characterizations (see \cite{kirchberg:invent},
\cite{haagerup-winslow}) but the results of this section show that the
difference between embeddability and hyperfiniteness is quite
delicate (at least for McDuff factors).

\begin{thm}
\label{thm:embeddableMcDuff}
If $M \subset B(L^2(M))$ is a McDuff factor with trace $\tau_M$ then
the following are equivalent:

\begin{enumerate}
\item $M$ is embeddable into $R^{\omega}$.

\item There exists a weakly dense C$^*$-subalgebra $A \subset M$ and a
sequence of normal, u.c.p. maps $\varphi_n : M \to M_{k(n)} ({\mathbb
C})$ such that

\begin{enumerate}
\item $\| \varphi_n (ab) - \varphi_n (a) \varphi_n (b) \|_2 \to 0$ and

\item $| \tau_{k(n)} \circ \varphi_n (a) - \tau_M (a) | \to 0$ for all
$a,b \in A$.
\end{enumerate}

\item There exists a weakly dense C$^*$-subalgebra $A \subset M$ and
finite rank projections $\{ P_n \}$  such that

\begin{enumerate}
\item $\frac{ \| [P_n, a] \|_{HS}}{\| P_n \|_{HS}} \to 0$  and 

\item $\frac{<aP_n, P_n>_{HS}}{<P_n,P_n>_{HS}} \to \tau_M (a)$ for all
$a \in A$.
\end{enumerate}

\item There exists a weakly dense C$^*$-subalgebra $A \subset M$ and a
state $\phi$ on $B(L^2(M))$ such that

\begin{enumerate}
\item $A \subset B(L^2(M))_{\phi}$ and 

\item $\phi|_A = \tau_M.$
\end{enumerate}

\item There exists a weakly dense C$^*$-subalgebra $A \subset M$ and
an embedding $\Phi : M \hookrightarrow R^{\omega}$ such that $\Phi|_A$
is completely positively liftable.

\item There exists a weakly dense C$^*$-subalgebra $A \subset M$ such
that $C^*(A, JAJ) \cong A \otimes A^{op}$.

\item $M$ has a weak expectation relative to a weakly dense 
C$^*$-subalgebra. (i.e. There exists a weakly dense 
C$^*$-subalgebra $A \subset M$ and 
a u.c.p.\ map $\Phi : B(L^2(M)) \to M$ such that $\Phi|_A = id_A$.

\item There exists a weakly dense operator system $X \subset M$ such 
that $X$ is injective.
\end{enumerate}
\end{thm}

\begin{proof} 
(1) $\Longrightarrow$ (2). Choose a trace $\gamma \in T(C^*({\mathbb
    F}_{\infty}))$ such that the weak closure of the GNS
    representation is isomorphic to $M$.  Since $M$ embeds into
    $R^{\omega}$ it follows that $\gamma$ is in the weak-$*$ closure of
    the set $T(C^*({\mathbb F}_{\infty}))_{FD}$ (traces whose GNS
    representations are finite dimensional).  This follows easily from
    the remark that unitaries in $R^{\omega}$ always lift to
    unitaries in $l^{\infty}(R)$ (see \cite[Lemma
    4.5]{kirchberg:invent}).  

From Theorem \ref{thm:basicconstruction} it follows that we can 
find a Popa algebra $A$ with a trace $\tau \in \TAafd$ such that the 
weak closure of the GNS representation of $A$ with respect to $\tau$ 
is isomorphic to $M$ ($\cong M \bar{\otimes} R$, since $M$ is McDuff).  
(2) now follows from the definition of $\TAafd$, Arveson's Extension 
Theorem and the remark that u.c.p. maps from a von Neumann algebra to 
a matrix algebra can always be approximated by normal, u.c.p. maps.

Using Theorem \ref{thm:mainthmsection3} it is easy to see that
statements (2) - (7) are equivalent.  Since, (5) obviously implies
(1), it follows that (1) - (7) are equivalent.  

(8) $\Longrightarrow$ (1) is a consequence of 
\cite[Theorem 1.4]{kirchberg:invent} together with the equivalence of 
($vi$) and ($iii$) in \cite[Proposition 1.3]{kirchberg:invent} (recall 
that $B(H)$ has the WEP). 

(7) $\Longrightarrow$ (8).  If $A \subset M$ is weakly dense and $\Phi 
: B(L^2(M)) \to M$ is a weak expectation relative to $A$ then 
\cite[Theorem 2.1]{blackadar:Proc.WEP} ensures that we can find an 
idempotent u.c.p.\ map $\Psi : B(L^2(M)) \to M$ such that $\Psi (a) = a$ 
for all $a \in A$.  The desired injective operator system is then 
$X = \Psi(B(L^2(M)))$ (since $\Psi \circ \Psi = \Psi$).
\end{proof}

In \cite{blackadar:WEP} Blackadar proved the existence of a
non-injective factor which has a weak expectation relative to a dense
subalgebra.  Other examples are exhibited by Kirchberg in the remark
after \cite[Corollary 3.5]{kirchberg:invent}.  However, in both of
these papers it is far from clear if finite examples can be
constructed using their techniques. Thus part (7) of Theorem
\ref{thm:embeddableMcDuff} seems to give the first examples of
non-hyperfinite II$_1$ factors with weak expectations and shows that
in fact many well known II$_1$ factors have weak expectations.

In order to illustrate just how delicate the results above are, we
remind the reader of the various characterizations of the hyperfinite
II$_1$ factor (most of which are due to Alain Connes).

\begin{thm}
\label{thm:R}
If $M \subset B(L^2(M))$ is a II$_1$ factor with trace $\tau_M$ then
the following are equivalent:

\begin{enumerate}
\item $M \cong R$.

\item There exists a weakly dense C$^*$-subalgebra $A \subset M$ and a
sequence of normal, u.c.p. maps $\varphi_n : M \to M_{k(n)} ({\mathbb
C})$ such that

\begin{enumerate}
\item $\| \varphi_n (ab) - \varphi_n (a) \varphi_n (b) \|_2 \to 0$ for
all $a,b \in A$ and

\item $| \tau_{k(n)} \circ \varphi_n (x) - \tau_M (x) | \to 0$ for all
$x \in M$.
\end{enumerate}

\item There exists a weakly dense C$^*$-subalgebra $A \subset M$ and
finite rank projections $\{ P_n \}$ such that

\begin{enumerate}
\item $\frac{ \| [P_n, a] \|_{HS}}{\| P_n \|_{HS}} \to 0$ for every $a
\in A$ and 

\item $\frac{<xP_n, P_n>_{HS}}{<P_n,P_n>_{HS}} \to \tau_M (x)$ for all
$x \in M$.
\end{enumerate}

\item There exists a weakly dense C$^*$-subalgebra $A \subset M$ and a
state $\phi$ on $B(L^2(M))$ such that

\begin{enumerate}
\item $A \subset B(L^2(M))_{\phi}$ and 

\item $\phi|_M = \tau_M.$
\end{enumerate}

\item There exists a completely positively liftable embedding $\Phi :
M \hookrightarrow R^{\omega}$.

\item $C^*(M, JMJ) \cong M \otimes M^{op}$.

\item There exists a u.c.p. map $\Phi : B(L^2(M)) \to M$ such that 
$\Phi|_M = id_M$ (i.e.\ $M$ is injective).
\end{enumerate}
\end{thm}

\begin{proof} (1) $\Longrightarrow$ (2) is obvious.  
(2) $\Longrightarrow$ (5) is contained in the argument  
in the proof of (1) $\Longrightarrow$ (2) from Theorem 
\ref{thm:mainthmUTAwafd}.  Since (5) is equivalent to hyperfiniteness 
for general finite von Neumann algebras, we see that (1), (2) and (5) 
are equivalent. 

The equivalence of (1), (6) and (7) is due to Connes (cf.\ 
\cite[Theorem 5.1]{connes:classification}).  This paper also contains 
his adaptation of Day's trick to deduce (3) from injectivity (7). 
Hence we are left to prove (3) $\Longrightarrow$ (4) and 
(4) $\Longrightarrow$ (1).

(3) $\Longrightarrow$ (4) is well known: simply take a cluster point
of the states $$T \mapsto \frac{<TP_n, P_n >_{HS}}{<P_n, P_n
>_{HS}}$$ on $B(L^2(M))$.  Finally, (4) $\Longrightarrow$ (1) also
follows from Connes' work since the density of $A$ in $M$, together
with the fact that the hypertrace takes the correct value on all of
$M$, implies that actually $M$ is contained in the centralizer of
$\varphi$ and hence is hyperfinite (cf.\ proof of \cite[Theorem
5.1]{connes:compactmetricspaces}).
\end{proof}

\section{Szeg\"{o}'s Limit Theorem for Self Adjoint Operators} 

Here we observe that Proposition \ref{thm:traciallytypeI} together
with Voiculescu's Theorem yield a general form of a classical theorem
of Szeg\"{o}.  We recommend \cite{arveson:survey},
\cite{arveson:numerical} and \cite{bedos:Szego} for nice discussions
of this problem, it's relation to numerical approximation of spectra
and overviews of the vast body of previous work on this problem. 
(We are particularly fond of the
introduction in \cite{bedos:Szego}.)

Let $T \in B(H)$ be a self adjoint operator, $A = C^* (T, I_H)$ be the
C$^*$-algebra generated by $T$ and the identity operator and $\varphi
\in S(A) = T(A)$ be an arbitrary (tracial) state.  Let $\mu_{\varphi}$
be the spectral distribution of $T$ with respect to $\varphi$; i.e.\
the unique probability measure on $\mathbb R$ such that
$$\int_{-\infty}^{+\infty} f(x) d\mu_{\varphi}(x) = \varphi(f(T)),$$
for all $f \in C_0 ({\mathbb R})$. 

If $\{ P_n \}$ is any sequence of finite rank projections then we will
regard the compressions $P_n T P_n$ as self adjoint operators on the
finite dimensional Hilbert space $P_n (H)$.  We will let $\mu_n$
denote the spectral distribution of $P_n T P_n$ with respect to the
unique trace on the matrix algebra in which $P_n T P_n$ sits.

\begin{thm}
\label{thm:Szego}
Let $T \in B(H)$ be a self adjoint operator and $A = C^* (T, I_H)$.
Then, there exists a sequence of finite rank projections $P_1 \leq P_2
\leq P_3 \ldots \in B(H\otimes_2 H)$ such that $P_n \to I_{H\otimes_2
H}$ in the strong operator topology and for every $\varphi \in \TA$
there exists a subsequence $\{ n_k \}$ such that $$\mu_{n_k} \to
\mu_{\varphi}$$ in the weak-$*$ topology, where $\mu_{n_k}$ denotes
the spectral distribution of $P_{n_k} T\otimes I_H P_{n_k}$ as above.
\end{thm}

\begin{proof} 
Since abelian C$^*$-algebras are type I, every state (i.e.\ trace) is
uniformly approximately finite dimensional.  Hence, since $A\otimes
I_H$ contains no non-zero compact operators, we can apply Proposition
\ref{thm:QDcase} to find a sequence of projections $P_1 \leq P_2 \leq
P_3 \ldots \in B(H\otimes_2 H)$ which asymptotically commute with
$T\otimes I_H$ in norm and recapture the entire state space of $A$.
The theorem now follows from \cite[Theorem 6]{bedos:Szego}.
\end{proof}

Note that if $A \cap {\mathcal K}(H) = \{0\}$ then the filtration can
be found on $H$ (i.e.\ one does not have to change the representation
in this case).  It should also be remarked that numerical analysts
will likely not find the above theorem to be of any use as the
construction of the projections $P_n$ requires knowledge of the
spectral distribution.

\section{Questions}

The following questions seem natural, in light of the present work. 

\begin{enumerate} 
\item Is every II$_1$ factor representation of a Popa algebra McDuff? 
While this seems unlikely, the inductive limit constructions 
used in the classification program tend to produce McDuff factors.

\item Is there an example such that $\TAwafd \neq \TAafd$ (or
$\UTAwafd \neq \UTAafd$)?  The only obvious obstruction is related to
quasidiagonality since the existence of a {\em faithful} trace in
$\TAafd$ implies quasidiagonality.  However, as we saw in Section 7,
$\TAwafd = \emptyset$ for any of the standard examples of stably finite,
non-QD C$^*$-algebras.  For example, can one construct a C$^*$-algebra
which has the WEP and a {\em faithful} trace but which is not QD? The 
most natural candidate seems to be the hyperfinite II$_1$ factor.

\item Can a free group factor or a II$_1$ factor with property T
contain a weakly dense, QD C$^*$-subalgebra? How about a Popa algebra?
Our constructions always give McDuff factors (so we have some place to
hide the Popa algebra).  Note that amenability can't be an obstruction
since $L(G_1 \times G_2) \cong L(G_1) \bar{\otimes} L(G_2)$ and hence
many non-amenable groups give McDuff factors and hence contain dense
QD subalgebras by Theorem \ref{thm:arbitraryMcDuff}.

\item Can one give estimates of the free entropy dimension of a finite
set of elements (not necessarily generators) in a Popa algebra which
is independent of the particular trace?  This is related to the
semicontinuity problem for free entropy dimension.

\item Can one prove a classification theorem for simple, nuclear, real
rank zero C$^*$-algebras which satisfy the UCT and such that $\TA =
\TAafd$ ($= \UTAafd$, by Theorem \ref{thm:locallyreflexive})? In
\cite[Theorem 3.3]{popa:simpleQD} Popa proves that such algebras have
an internal finite dimensional approximation property which should be
of use.  Presumably the role of nuclearity needs to be clarified as
Popa never assumes nuclearity in \cite{popa:simpleQD}.

\item Is every (nuclear) Popa algebra with real rank zero, unperforated
K-theory, Riesz decomposition property and unique trace (or, perhaps,
finitely many extreme traces) necessarily tracially AF?  

\item Can an infinite, simple, discrete group with Kazdan's property T
be embed into the unitary group of an $R^{\omega}$-embeddable McDuff
factor? (Compare with \cite{robertson} where it is shown that no such
embedding exists into the unitary group of $L({\mathbb
F}_n)\bar{\otimes} R$ or, more generally, $L(\Gamma)$ for any a-T-amenable 
discrete group $\Gamma$.) 

\item Let $\Gamma$ be a discrete group such that the group von Neumann
algebra of $\Gamma$ contains a weakly dense operator system which is
injective.  Does it follow that $\Gamma$ is an exact group? Perhaps
just uniformly embeddable into Hilbert space?  Since these notions
pass to subgroups, this would imply that every residually finite group
(and every other group which embeds into the unitary group of
$R^{\omega}$) is exact (or uniformly embeddable). It would also show
that non-embeddable groups, whose existence has been asserted by
Gromov, give counterexamples to Connes' embedding problem.
\end{enumerate}

\section{Appendix: Is $R$ Quasidiagonal?} 

Here we discuss a basic problem in operator algebra theory about which
very little is known.  Namely, whether or not the hyperfinite II$_1$
factor is QD.  One has to be careful as we are now thinking of $R$ as
a C$^*$-algebra and as such it is no longer separable (in norm).  This
has led to some confusion as the `classical' definition of
quasidiagonality is not the correct notion for non-separable
C$^*$-algebras. Voiculescu gives the correct local definition in the 
non-separable case in \cite{dvv:QDsurvey}. (See also the appendix to 
\cite{brown:QDsurvey} for a detailed treatment of the non-separable case.)

We are not aware of a reference for the following fact, however it is
known to a number of experts.  We thank George Elliott for showing us
a very nice proof.  We have only modified the last few lines of
Elliott's argument so that we can deduce a slightly stronger
statement.

\begin{lem}
Let $R$ act on $L^2(R)$ via the GNS construction.  There is no
 {\em sequence} of nonzero, finite rank projections $P_1, P_2,
 \ldots$ such that $\| [x,P_n] \| \to 0$ for all $x \in R$. 
\end{lem}
 
\begin{proof}  The proof goes by contradiction.  So let $P_1, P_2,
\ldots$ be finite rank projections such that $\| [x,P_n] \| \to 0$ for
all $x \in R$.  Put $K = \oplus_{n \in \mathbb N} L^2(R) = L^2(R)
\otimes_2 l^2({\mathbb N})$ and $P = \oplus_{n \in \mathbb N} P_n$.
Then $(x\otimes 1)P - P(x\otimes 1)$ is a compact operator for every
$x \in R$.  Hence, down in the Calkin algebra $P$ will land in the
commutant of $R\otimes 1$.  But then by a theorem of Johnson and
Parrott (see the remarks after \cite[Lemma 3.3]{johnson-parret}) it
follows that $P$ is a compact perturbation of an element in the
commutant of $R\otimes 1$. That is, there exists an infinite matrix $T
= (T_{i,j})_{i,j \in {\mathbb N}}$ such that each $T_{i,j} \in
R^{\prime} \subset B(L^2(R))$ and $P - T$ is compact on $K$.  In
particular, this implies that $\| P_n - T_{n,n} \| \to 0$.  Thus $\|
T_{n,n} \| \to 1$ and down in the Calkin algebra the norm of $T_{n,n}$
is tending to zero.  However this is a contradiction since the
commutant of $R$ is a II$_1$ factor (isomorphic to $R$) and hence a
simple C$^*$-algebra.  Thus the mapping to the Calkin algebra is
isometric.
\end{proof}

Note that the proof above never used the fact that the $P_n$'s are
projections and hence also holds for {\em sequences} of finite rank
operators whose norms are tending to one.  However, since we can
always construct a quasicentral {\em net} of finite rank operators for
$R$ we are left to conclude that the lemma above has more to do with
sequences versus nets (i.e.\ separable versus non-separable Hilbert
spaces) than it does with quasidiagonality.

The only other obvious strategy for proving nonquasidiagonality of $R$
is to embed a non-QD C$^*$-algebra into it.  However we already
observed in Section 7 that none of the known examples of non-QD
C$^*$-algebras can be embed into $R$. 

The tracial invariants studied in this paper lead to another approach
for proving that $R$ is not QD. They also show that an affirmative
answer to this question would have some remarkable consequences.

\begin{prop}
$R$ is QD if and only if for every (separable) C$^*$-algebra $A$ we
have $\TAafd \supset \UTAwafd$.
\end{prop}
 
\begin{proof}
We begin with the necessity. Let $A$ be arbitrary.  It suffices to
show that the extreme points of $\UTAwafd$ belong to $\TAafd$.
However, every extreme point of $\UTAwafd$ is also an extreme point of
$\TA$, since $\UTAwafd$ is a face, and hence gives $R$ in the GNS
representation.  Thus, if we assume that $R$ is QD then it's unique
trace must belong to $T(R)_{\rm AFD}$ and this completes the proof. 

For the sufficiency, we first point out that $R$ is QD if and only if
all of it's separable C$^*$-subalgebras are QD.  So let $A \subset R$
be an arbitrary separable, unital subalgebra.  Let $\tau \in \TA$ be
the restriction of the unique trace on $R$ to $A$.  Clearly $\tau$ is
faithful and belongs to $\UTAwafd$.  Hence it also belongs to
$\TAafd$.  This implies that $A$ is QD by Voiculescu's abstract
characterization of quasidiagonality (cf.\ \cite[Theorem
1]{dvv:QDhomotopy}).
\end{proof}

Applying the proposition above and Theorem \ref{thm:locallyreflexive}
we immediately get the following corollary.

\begin{cor}
If $R$ is QD then for every locally reflexive C$^*$-algebra $A$ we
have $$\TAwafd = \UTAwafd = \TAafd = \UTAafd.$$
\end{cor}

\bibliographystyle{amsplain}

\providecommand{\bysame}{\leavevmode\hbox to3em{\hrulefill}\thinspace}

\end{document}